\documentclass[preprint,
               11pt,
               authoryear]{elsarticle}

\usepackage[inline]{trackchanges} 
\usepackage{soul}
\usepackage{subfiles}
\usepackage{subcaption}
\usepackage{commath}
\usepackage{comment}
\usepackage{enumitem}
\usepackage{nicefrac}
\usepackage{booktabs}
\usepackage{xfp}
\usepackage{lineno}
\usepackage{amssymb}
\usepackage{url}
\usepackage{amsthm}
\usepackage{tikz}
\usetikzlibrary{
    shapes, 
    calc,
    decorations.pathreplacing, 
    calligraphy, 
    positioning,
    backgrounds
}

\addeditor{Jeremy}
\addeditor{Giacomo}
\addeditor{Mark}
\addeditor{Bob}
\addeditor{Darren}

\soulregister\citep7
\soulregister\ref7
\soulregister\citet7
\soulregister\eqref7

\renewcommand{\phi}{\varphi}
\renewcommand{\O}{\mathcal{O}}
\renewcommand{\u}{\mathbf{u}}
\newcommand{\h}{\mathbf{h}}

\newcommand{\numax}{\nu^\text{max}}

\newcommand{\F}{\mathcal{F}}
\newcommand{\C}{\mathcal{C}}
\newcommand{\ifone}{\text{IF-1}}
\newcommand{\iftwo}{\text{IF-2}}
\newcommand{\interior}{\text{int}}
\newcommand{\Ord}[1]{\O\left((\Delta t)^{#1}\right)}

\newtheorem{theorem}{Theorem}
\newdefinition{remark}{Remark}

\definecolor{myred}{HTML}{ee4035}
\definecolor{myorange}{HTML}{f37736}
\definecolor{myyellow}{HTML}{fdf498}
\definecolor{mygreen}{HTML}{7bc043}
\definecolor{myblue}{HTML}{0392cf}

\definecolor{coarse}{HTML}{ff0000}
\definecolor{coarseoutline}{HTML}{8b0000}

\definecolor{int2}{HTML}{ffd700}
\definecolor{int2outline}{HTML}{ffa500}

\definecolor{int1}{HTML}{ffc0cb}
\definecolor{int1outline}{HTML}{ee82ee}

\definecolor{fineish}{HTML}{e0ffff}
\definecolor{fineishoutline}{HTML}{add8e6}

\definecolor{fine}{HTML}{40e0d0}
\definecolor{fineoutline}{HTML}{0000ff}

\journal{Journal of Computational Physics}

\begin{document}


\begin{frontmatter}

    \title{Local Time-Stepping for the Shallow Water Equations using CFL Optimized Forward-Backward Runge-Kutta Schemes}

    \author[osu,ccs2]{Jeremy R. Lilly}
    \author[ccs2]{Giacomo Capodaglio}
    \author[t3]{Darren Engwirda}
    \author[osu]{Robert L. Higdon}
    \author[ccs2]{Mark R. Petersen}

    \affiliation[osu]{organization={Department of Mathematics},
            addressline={Oregon State University}, 
            city={Corvallis},
            state={OR},
            postcode={97331}, 
            country={USA},
    }
    \affiliation[ccs2]{organization={Computational Physics and Methods Group},
            addressline={Los Alamos National Laboratory}, 
            city={Los Alamos},
            state={NM},
            postcode={87545}, 
            country={USA},
    }
    \affiliation[t3]{organization={Fluid Dynamics and Solid Mechanics Group},
            addressline={Los Alamos National Laboratory}, 
            city={Los Alamos},
            state={NM},
            postcode={87545}, 
            country={USA},
    }

    \begin{abstract}
        The Courant–Friedrichs–Lewy (CFL) condition is a well known, necessary condition for the stability of explicit time-stepping schemes that effectively places a limit on the size of the largest admittable time-step for a given problem.
        We formulate and present a new local time-stepping (LTS) scheme optimized, in the CFL sense, for the shallow water equations (SWEs).
        This new scheme, called FB-LTS, is based on the CFL optimized forward-backward Runge-Kutta schemes from \citet{lilly2023_cfl}.
        We show that FB-LTS maintains exact conservation of mass and absolute vorticity when applied to the TRiSK spatial discretization \citep{ringler2010}, and provide numerical experiments showing that it retains the temporal order of the scheme on which it is based (second order).
        In terms of computational performance, we show that when applied to a real-world test case on a highly-variable resolution mesh, the MPAS-Ocean implementation of FB-LTS is up to 10 times faster than the classical four-stage, fourth-order Runge-Kutta method (RK4), and 2.3 times faster than an existing strong stability preserving Runge-Kutta based LTS scheme (LTS3).
        Despite this significant increase in efficiency, the solutions produced by FB-LTS are qualitatively equivalent to those produced by both RK4 and LTS3.
    \end{abstract}


\end{frontmatter}




\section{Introduction}
\label{sec:introduction}

The computational performance of explicit time-stepping schemes is often limited by the so-called Courant–Friedrichs–Lewy (CFL) condition.
Given a system of partial differential equations (PDEs) with a finite speed of propagation and a spatial discretization, the CFL condition bounds the time-step above in terms of the size of the spatial cells and a quantity related to the speed of the dynamics of the problem.
Formally, the CFL condition states that it is necessary for stability that
\begin{linenomath}
\begin{equation}
    \label{eqn:cfl}
    \nu = c\frac{\Delta t}{\Delta x} \leq \numax \,,
\end{equation}
\end{linenomath}
where \( c \) is a characteristic speed, such as that of a gravity wave, \( \nu \) is referred to as the Courant number, and \( \numax \) is the maximum admittable Courant number, which depends on the model problem itself and the chosen time and space discretizations.
For many applications, particularly those that require high performance computing (HPC) resources, the size of the desired spatial discretization \( \Delta x \) and the speed \( c \) can vary wildly across the computational domain.
In these cases, traditional explicit global time-stepping schemes are forced to satisfy a single, global CFL condition, perhaps enforced by only a small portion of the domain.
Local time-stepping (LTS) methods provide an answer to this problem by allowing a scheme to take time-steps depending on local values of \( c \) and \( \Delta x \), satisfying a local CFL condition rather than a global one.
In practice, this means that a domain can be divided into \textit{coarse} regions where large time-steps are used, and \textit{fine} regions, where smaller time-steps are used.
In contrast, a global scheme requires the small time-step native to a fine region to be used everywhere, resulting in the need for more evaluations of right-hand-side terms to advance forward in time, increasing the overall computational burden.

In this work, we introduce a new LTS scheme for the shallow water equations (SWEs) that has been optimized in the CFL sense. 
This scheme, called FB-LTS, is based on the three-stage, second-order, forward-backward Runge-Kutta scheme FB-RK(3,2) developed by \citet{lilly2023_cfl}.
FB-RK(3,2) is an explicit time-stepping scheme designed for coupled systems of PDEs. 
In the context of the SWEs, the scheme uses a forward-backward (FB) average of available thickness data to advance the momentum equation at each Runge-Kutta stage.
The weights of these FB averages have been optimized to produce a scheme that has maximal \( \numax \) when applied to the SWEs.
It was shown in \citet{lilly2023_cfl} that FB-RK(3,2) outperforms a popular three-stage, third-order strong stability preserving Runge-Kutta scheme (SSPRK3) in admittable time step by factors roughly between 1.6 and 2.2, making the scheme approximately twice as computationally efﬁcient with little to no effect on solution quality.
The new FB-LTS scheme introduced here takes advantage of the CFL performance of FB-RK(3,2) and combines it with the benefits of a local time-stepping scheme.

The particular algorithm by which regions of the domain using different time steps communicate,
which we refer to here as the LTS framework, was originally developed by \citet{hoang2019} for use with a TRiSK spatial discretization \citep{ringler2010}. 
TRiSK is a finite volume-type spatial discretization made for unstructured, variable-resolution polygonal grids, and is the discretization used in the Model for Prediction Across Scales-Ocean (MPAS-O) \citep{ringler2010, ringler2013, petersen2019}.
The application of this LTS framework in the context of a TRiSK spatial discretization is a major topic of this work.
The scheme presented in \citet{hoang2019} is based on SSPRK3, and is referred to here as LTS3.
Using FB-RK(3,2) as opposed to SSPRK3, our FB-LTS scheme is able to outperform LTS3 in terms of the size of the admittable time-step by factors up to 2.3.

The long-term goal of FB-LTS is to increase the computational efficiency of climate-scale models of the ocean and atmosphere running on highly variable resolution meshes, 
with a particular focus on the Energy Exascale Earth System Model (E3SM) being developed by the U.S. Department of Energy \citep{golaz2022}.
In this work, we implement FB-LTS, along with a certain operator splitting, for single-layer configurations in MPAS-Ocean, the ocean component of E3SM.
A single-layer ocean is modeled by the SWEs, and serves as the starting point for our eventual goal to use FB-LTS in multi-layer, climate-scale models.

This paper is structured as follows.
We recall the formulation of FB-RK(3,2) presented in \citet{lilly2023_cfl} for completeness, then present the FB-LTS scheme.
Next, we show that the scheme exactly conserves mass and absolute vorticity in the context of a TRiSK spatial discretization.
Then, we discuss the details of the implementation of the scheme in MPAS-Ocean, including a discussion of an operator splitting approach that we have adopted within the SWEs.
Finally, we perform a number of numerical experiments that demonstrate the computational efficiency of FB-LTS as compared to LTS3 in MPAS-Ocean.
These experiments model the storm surge in Delaware Bay caused by hurricane Sandy in 2012, and are an evolution of the simulations explored in \cite{lilly2023_sandy}.


\section{Local Time-Stepping Schemes with FB-RK(3,2)}
\label{sec:fblts}

We begin by recalling the FB-RK(3,2) scheme introduced by \citet{lilly2023_cfl} for completeness, then introduce a new LTS scheme based on this global scheme, called FB-LTS.
The primary goal of FB-LTS is to solve the shallow water equations (SWEs) efficiently in the CFL sense, i.e. taking time-steps as large as possible.
To facilitate discussion of our methods and the SWEs, we introduce the nonlinear SWEs on a rotating sphere, given by
\begin{linenomath}
\begin{equation}
\begin{aligned}
    \pd{\u}{t} + \left( \nabla \times \u + f\mathbf{k} \right) \times \u &= -\nabla\frac{\abs{\u}^2}{2} - g\nabla (h + z_b) \\
    \pd{h}{t} + \nabla \cdot \left(h\u\right) &= 0 \,,
\end{aligned} \label{eqn:nonlinear_swe}
\end{equation}
\end{linenomath}
where \( \u(x, y, t) = \left(u(x, y, t),\, v(x, y, t)\right)\) is the horizontal fluid velocity,
\( x \) and \( y \) are the spatial coordinates,
\( t \) is the time coordinate,
\( f \) is the Coriolis parameter,
\( \mathbf{k} \) is the local vertical unit vector,
\( g \) is the gravitational constant,
\( h \) is the fluid thickness,
and \(z_{b}\) is the height of the bottom topography.
In Section \ref{subsec:hurricane_model}, we will introduce a similar shallow water model of particular interest that will showcase the performance of FB-LTS in a real-world test case.
Throughout this work, we often refer to a given equation for the evolution of \( \u \) as the momentum equation, and a given equation for the evolution of \( h \) as the thickness or mass equation.


\subsection{FB-RK(3,2)}
\label{subsec:fbrk32}

The time-stepping scheme presented here is an extension of the three-stage, second-order Runge-Kutta time-stepping scheme RK(3,2) from \citet{wicker2002} which is used to solve a SWE-like system in MPAS-Atmosphere.
This extension of RK(3,2) allows the use of the most recently obtained data for the layer thickness to update the momentum data within each Runge-Kutta stage.
This is done by taking a weighted average of layer thickness data at the old time level \( t^n \) and the most recent RK stage, then applying this to the momentum equation.

Consider a general system of ODEs in independent variables \( u = u(t) \) and \( h = h(t) \) of the form
\begin{linenomath}
\begin{equation}
\begin{aligned}
    \od{u}{t} &= \Phi\left( u, h \right) \\
    \od{h}{t} &= \Psi\left( u, h \right) \,,
\end{aligned} \label{eqn:odesys}
\end{equation}
\end{linenomath}
where \( t \) is the time coordinate.
As discussed above, in the context of the SWEs \( u \) is the fluid velocity and \( h \) is the ocean layer thickness.
The right-hand-side operators \( \Phi \) and \( \Psi \) are refereed to as the momentum and thickness (or mass) tendencies respectively.
Let \( u^n \approx u(t^n)\) and \( h^n \approx h(t^n) \) be the numerical approximations to \( u \) and \( h \) at time \( t = t^n \).
Let \( \Delta t \) be a time-step such that \( t^{n+1} = t^n + \Delta t \).
Also, let \( t^{n+\nicefrac{1}{m}} = t^n + \frac{\Delta t}{m} \) for any positive integer \( m \).
Then, FB-RK(3,2) is given by
\begin{linenomath}
\begin{subequations}
\label{eqn:fbrk32}
\begin{align}
\begin{split}
    \bar{h}^{n+\nicefrac{1}{3}} &= h^n + \frac{\Delta t}{3} \Psi\left( u^n, h^n \right) \\
    \bar{u}^{n+\nicefrac{1}{3}} &= u^n + \frac{\Delta t}{3} \Phi\left( u^n, h^* \right) \\
    h^* &= \beta_1 \bar{h}^{n+\nicefrac{1}{3}} + (1-\beta_1) h^n
\end{split} \label{subeqn:fbrk32_stage1} \\
\nonumber \\
\begin{split}
    \bar{h}^{n+\nicefrac{1}{2}} &= h^n + \frac{\Delta t}{2} \Psi\left( \bar{u}^{n+\nicefrac{1}{3}}, \bar{h}^{n+\nicefrac{1}{3}} \right) \\
    \bar{u}^{n+\nicefrac{1}{2}} &= u^n + \frac{\Delta t}{2} \Phi\left( \bar{u}^{n+\nicefrac{1}{3}}, h^{**} \right) \\
    h^{**} &= \beta_2 \bar{h}^{n+\nicefrac{1}{2}} + (1-\beta_2) h^n
\end{split} \label{subeqn:fbrk32_stage2} \\
\nonumber \\
\begin{split}
    h^{n+1} &= h^n + \Delta t \Psi\left( \bar{u}^{n+\nicefrac{1}{2}}, \bar{h}^{n+\nicefrac{1}{2}} \right) \\
    u^{n+1} &= u^n + \Delta t \Phi\left( \bar{u}^{n+\nicefrac{1}{2}}, h^{***} \right) \\
    h^{***} &= \beta_3 h^{n+1} + (1-2\beta_3) \bar{h}^{n+\nicefrac{1}{2}} + \beta_3 h^n \,.
\end{split} \label{subeqn:fbrk32_stage3}
\end{align}
\end{subequations}
\end{linenomath}
The weights \( \beta_1 \), \( \beta_2 \), and \( \beta_3 \) are called the forward-backward (FB) weights.
These FB-weights can be chosen so as to optimize the allowable time-step in the SWEs; it was shown in \cite{lilly2023_cfl} that taking \( (\beta_1,\, \beta_2,\, \beta_3) = (0.531,\, 0.531,\, 0.313) \) increases the admittable time-step versus RK(3,2) between factors of 1.6 and 2.2 in a number of nonlinear test cases while maintaining second-order accuracy.


\subsection{FB-LTS}
\label{subsec:fblts}

\begin{figure}
    \centering
    
    \def\dotscale{1.3}
    \def\mydot{\scalebox{\dotscale}{\( \mathcolor{black}{\bullet} \)}}
    \def\mysquare{\scalebox{\dotscale}{{\tiny \( \mathcolor{black}{\blacksquare} \)}}}
    \def\mytriangle{\scalebox{\dotscale}{\raisebox{1px}{\small\( \mathcolor{black}{\blacktriangle} \)}}}

    \tikzstyle{tn}=[font=\Large,align=center]
    \tikzstyle{ar}=[->, dotted, very thick]
    \tikzstyle{background grid}=[draw, black!50, step=1cm]
    \begin{tikzpicture}[scale=3]
        \coordinate (primal_center) at (0,0);
        \coordinate (dual_center) at (30:1);
        \coordinate (edge_center) at (60:0.866);
        \coordinate (mid_norm) at (60:0.5);

        \node[tn] at ($(edge_center)+(110:0.15cm)$) {\( \mathbf{x}_e \)};
    
        \node[tn, below left=-1pt of primal_center] {\( \mathbf{x}_i \)};
        \foreach \x in {30,90,...,330}{
            \draw[ultra thick] (\x:1) -- (\x+60:1);
        }
        \draw[ultra thick] (dual_center) -- ( $(dual_center)+(dual_center)$ );

        \node[tn, below right=1pt of dual_center] {\( \mathbf{x}_v \)};
        \foreach \x in {-150,-30,90}{
            \draw[ultra thick, dashed] ($(\x:1cm)+(dual_center)$) -- ($(\x+120:1cm)+(dual_center)$);
            \node at ($(\x:1cm)+(dual_center)$) {\mydot};
        }

        \draw[->, blue, line width=2pt] (edge_center) -- (mid_norm);
        \draw[rotate around={90:(edge_center)}, ->, red, line width=2pt] (edge_center) -- (60:0.5);
        \node[tn, below right=1pt of mid_norm] {\( \mathbf{n}_e \)};
        \node[tn, above left=-2pt of mid_norm, rotate around={90:(edge_center)}] {\rotatebox{-90}{\( \mathbf{t}_e \)}};

        \foreach \x in {30,90,...,330}{
            \node at (\x+60:1)  {\mytriangle};
            \node[rotate=\x-150] at (\x+30:0.866)  {\mysquare};
        }
        \node[rotate=30] at ($(30:0.5cm)+(dual_center)$)  {\mysquare};

        \node[tn] (Pi_label) at (120:1.5cm) {Primal cell \( P_i \)};
        \draw[ar] (Pi_label) edge[bend right=10] (150:0.5);

        \node[tn] (Dv_label) at (-15:1.75cm) {Dual cell \( D_v \)};
        \draw[ar] (Dv_label) edge[bend right=10] (15:1.25);

        \node[align=left, font=\large, text width=5.5cm, anchor=west] (legend) at (45:2) {
            \( \mathbf{x}_i = \) thickness points \newline
            \( \mathbf{x}_e = \) normal velocity points \newline
            \( \mathbf{x}_v = \) potential vorticity points
        };
    \end{tikzpicture}
    \caption{
        Example TRiSK grid from a Voronoi tessellation, where the primal cells are hexagons and the dual cells are triangles centered at primal cell vertices.
        This is the type of spatial discretzation used by MPAS-Ocean.
        The vector \( \mathbf{n}_e \) is normal to cell edge \( e \) in a fixed, arbitrary direction.
        Later, in Section \ref{subsec:conservation}, we define a quantity \( n_{e,i} \) as either 1 or -1 so that \( n_{e,i} \mathbf{n}_e \) is the outward unit normal vector to cell \( i \) at edge  \( e \).  
        Then, \( \mathbf{t}_e = \mathbf{k} \times \mathbf{n}_e \).
    }
    \label{fig:trisk}
\end{figure}
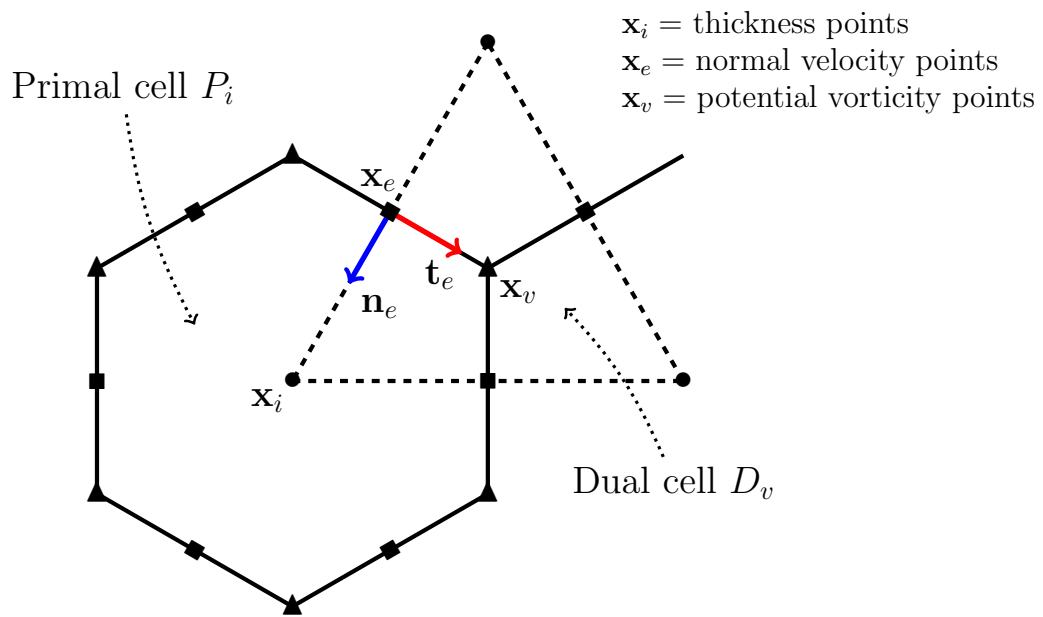

Here, we introduce FB-LTS in the context of a TRiSK spatial discretization \citep{ringler2010}, which is a finite volume-type spatial discretization made for unstructured, variable-resolution polygonal grids, and is the discretization used in MPAS-Ocean.
TRiSK employs C-grid-type discretization \citep{arakawa1977} wherein the mass variable is computed on cell centers and the normal component of velocity is computed on cell edges.
In MPAS-Ocean, these are Voronoi grids \citep{ju2011, okabe2017} consisting primarily of hexagons as the primal mesh, with a dual mesh consisting of triangles (Figure \ref{fig:trisk}).

Given some computational domain, let \( \Omega_P \) be the set of indices for primal cells (hereafter referred to as just cells, dual cells will be refereed to specifically as dual cells) and  \( \Omega_E \) be the set of indices for cell edges.
We decompose the computational domain into two regions, a coarse region, which will be advanced with the coarse time-step \( \Delta t \), and a fine region, which will be advanced with the fine time-step \( \frac{\Delta t}{M} \) for some positive integer \( M \).
Note that the label of fine or coarse does not necessarily reference the size of the spatial discretization, but rather which time-step is used to advance it.
The region made up of fine cells is called the fine region while the rest of the mesh is called the coarse region.
While the sets \( \Omega_P \) and \( \Omega_E \) are formally sets of indices, throughout the text we often refer to them directly as sets of cells and edges respectively for readability.

Following the notation of \citet{hoang2019}, let \( \F_P \) be the set of cells in the fine region, and \( \C_P \) be the set of cells in the coarse region such that \( \F_P \cup \C_P = \Omega_P \).
In the fine region, we define subsets \( \F_P^\ell \subseteq \F_P \) to be the so-called, interface adjacent fine cells.
Let \( r \) be the radius of the discrete tendency operators for \( \u \) and \( \h \), and define \( \F_P^\ell \) for \( \ell = 1, \cdots, 5 \) so that \( \F_P^\ell \) contains \( \ell r \) layers of cells in \( \F_P \) neighboring \( \C_P \).
For example, for FB-LTS as applied to a TRiSK spatial discretization, we have \( r = 2 \), so \( \F_P^5 \) contains 10 layers of cells bordering the interface one region (Figure \ref{fig:fblts_regions}).
These \( \F_P^\ell \) subsets are not disjoint with one another, rather \( \F_P^1 \subseteq \cdots \subseteq \F_P^5 \).

In the coarse region, define disjoint subsets \( \C_P^\ifone \subseteq \C_P \), \( \C_P^\iftwo \subseteq \C_P \), and \( \C_P^\interior \subseteq \C_P \) such that \( \C_P^\ifone \cup \C_P^\iftwo \cup \C_P^\interior = \C_P \).
Call \( \C_P^\ifone \) the set of interface one cells, \( \C_P^\iftwo \) the set of interface two cells, and \( \C_P^\interior \) the set of interior coarse cells; all these cells advance with the coarse time-step.
These collections of cells are distributed in the computational domain such that only \( \C_P^\ifone \) cells border \( \F_P \) cells, only \( \C_P^\iftwo \) cells border \( \C_P^\ifone \) cells, and only \( \C_P^\interior \) cells border \( \C_P^\iftwo \) cells.

The sets \( \F_E \), \(\F_E^\ell\), \( \C_E^\ifone \), \( \C_E^\iftwo \), \( \C_E^\interior \), and \( \C_E \) give the corresponding sets of cell edges for all the sets of cells described above.
An edge shared by a cell from \( \F_P \) and a cell from \( \C_P^\ifone \) belongs to \( \F_E \), an edge shared by a cell from \( \C_P^\ifone \) and a cell from \( \C_P^\iftwo \) belongs to \( \C_E^\ifone \), and an edge shared by a cell from \( \C_P^\iftwo \) and a cell from \( \C_P^\interior \) belongs to \( \C_E^\iftwo \).
In plain language, an edge in dispute between two cells of different regions belongs to the region closest to the fine region.
This domain decomposition is visualized in Figure \ref{fig:lts_regions}; we refer to these regions collectively as the LTS regions.
Often, we use the notation \( \C \) or \( \F \) without a subscript \( P \) or \( E \) to refer to the whole of the corresponding region, including both cells and edges.

Finally, we assume that the LTS regions are configured in such a way that there are enough layers of \( \C_P^\ifone \) and \( \C_P^\iftwo \) cells so that the operator stencils of the tendencies do not contain cells more than one region away.
For example, the operator stencil on fine cells can only contain fine cells and interface one cells, and an operator stencil on interface one cells can only contains fine, interface one, and interface two cells, Figure \ref{fig:fblts_lts3_regions} shows a practical example of this, where the radius of the operator stencil is \( r = 2 \).

\begin{figure}
    \centering
    \tikzstyle{background grid}=[draw, black!50, step=1cm]
    \tikzstyle{tn}=[very thick, rounded corners=0.25cm, font=\footnotesize, anchor=west, text width=4cm, align=center]
    \tikzstyle{ar}=[very thick, dotted, ->]
    \begin{tikzpicture}
        \node[opacity=1] (picture) at (0,0) {\includegraphics[width=0.6\textwidth]{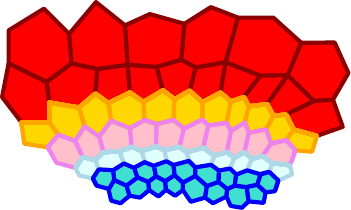}};

        \node[tn, fill=coarse, draw=coarseoutline] (coarse) at (5.5, 3) {\( \C^\interior \): coarse region, uses coarse time-step \( \Delta t \)};
        \node[tn, fill=int2, draw=int2outline] (int2) at (5.5, 1.5) {\( \C^\iftwo \): interface two region, uses coarse time-step \( \Delta t \)};
        \node[tn, fill=int1, draw=int1outline] (int1) at (5.5, 0) {\( \C^\ifone \): interface one region, uses coarse time-step \( \Delta t \)};
        \node[tn, fill=fineish, draw=fineishoutline] (fineish) at (5.5,-1.5) {\( \F^\ell \): interface adjacent fine region, uses fine time-step \( \nicefrac{\Delta t}{M} \)};
        \node[tn, fill=fine, draw=fineoutline] (fine) at (5.5,-3) {\( \F \): fine region, uses fine time-step \( \nicefrac{\Delta t}{M} \)};

        \node (coarse_cell) at (-0.6,1.7) {};
        \node (int2_cell) at (0.35,-0.15) {};
        \node (int1_cell) at (1.2,-0.95) {};
        \node (fineish_cell) at (2.2,-1.4) {};
        \node (fine_cell) at (1.6,-1.9) {};

        \path[ar]  (coarse_cell) edge[bend left=20] (coarse);
        \path[ar]  (int2_cell) edge[bend left=15] (int2);
        \path[ar]  (int1_cell) edge[bend left=10] (int1);
        \path[ar]  (fineish_cell) edge[bend right=5] (fineish);
        \path[ar]  (fine_cell) edge[bend right=10] (fine);
    \end{tikzpicture}   
    
    \caption{
        An example mesh with cells and edges labeled for LTS. Blue cells and edges belong to \( \F \), light blue cells and edges belong to \( \F^\ell \subseteq \F \), pink cells and edges belong to \( \C^\ifone \), yellow cells and edges belong to \( \C^\iftwo \), and red cells and edges belong to \( \C^\interior \).
        Note that in practice, one often needs more layers of light blue, pink, and yellow cells; see Figure \ref{fig:fblts_lts3_regions} for a practical example.
    }
    \label{fig:lts_regions}
\end{figure}

Now, consider the following system of PDE that has been discretized in space
\begin{linenomath}
\begin{equation}
\begin{aligned}
    \pd{u_e}{t} &= \Phi_e\left( \u, \h \right) \\
    \pd{h_i}{t} &= \Psi_i\left( \u, \h \right) \,,
\end{aligned} \label{eqn:semi_discrete_sys}
\end{equation}
\end{linenomath}
where \( \u = (u_e)_{e \in \Omega_E} \) and \( \h = (h_i)_{i \in \Omega_P} \).
Under a TRiSK spatial discretization, \( h_i \) is computed at primal cell centers, and \( u_e \) is the normal component of velocity computed at primal cell edges (Figure \ref{fig:trisk}).

Set a time-step \( \Delta t \), and let \( M \) be some positive integer.
Let \( t^{n+\nicefrac{1}{m}} = t^n + \frac{\Delta t}{m} \) for any positive integer \( m \), and \( t^{n,k} = t^n + k\frac{\Delta t}{M} \), and \( t^{n,k+\nicefrac{1}{m}} = t^n + \left(k+\nicefrac{1}{m} \right)\frac{\Delta t}{M} \) (Figure \ref{fig:time_levels}).
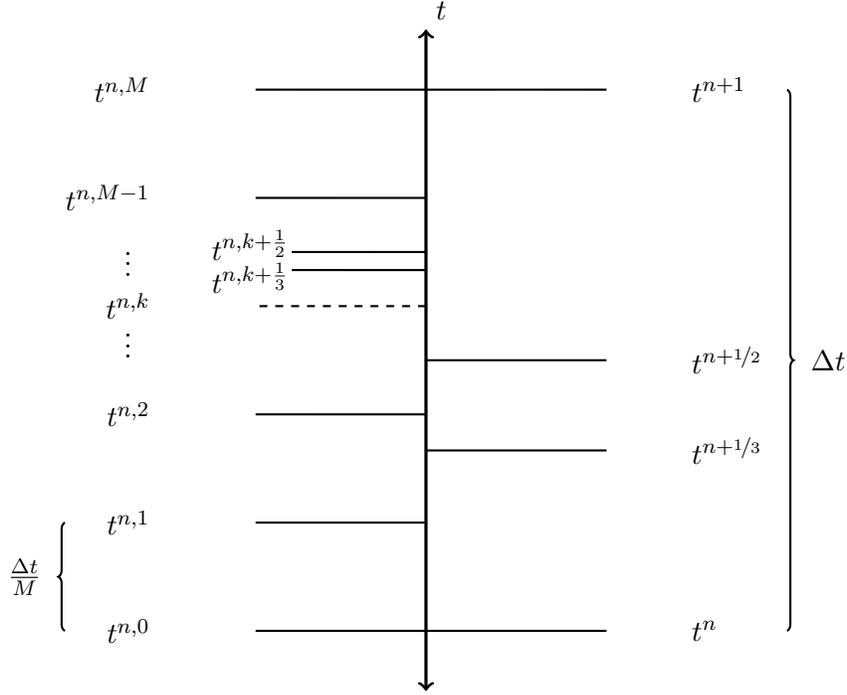
\begin{figure}
    \centering
    \def\width{3}
    \def\height{9}
    \def\numfine{5}
    \def\tnk{\inteval{\numfine-2}}
    \begin{tikzpicture}[scale=0.8]
        \coordinate (top) at (0,\height+1);
        \node[above right=0pt of top] {\( t \)};
        \draw[very thick, <->] (0,-1) -- (0,\height+1);

        \coordinate (coarse_start) at (\width,0);
        \draw[thick] (0,0) -- (coarse_start);
        \node[right=of coarse_start] (tn) {\( t^n \)};

        \coordinate (coarse_first) at (\width,\height/3);
        \draw[thick] (0,\height/3) -- (coarse_first);
        \node[right=of coarse_first] {\( t^{n+\nicefrac{1}{3}} \)};

        \coordinate (coarse_second) at (\width,\height/2);
        \draw[thick] (0,\height/2) -- (coarse_second);
        \node[right=of coarse_second] {\( t^{n+\nicefrac{1}{2}} \)};

        \coordinate (coarse_new) at (\width,\height);
        \draw[thick] (0,\height) -- (coarse_new);
        \node[right=of coarse_new] (tnp1) {\( t^{n+1} \)};
        
        \draw[decoration={brace,mirror,raise=0pt},decorate,thick] (\width+3,0) -- node[right=5pt] {\( \Delta t \)} (\width+3,\height);
        
        \foreach \i in {0,...,\numfine}{
            \node (fine_\i) at (-\width,\height*\i/\numfine) {};
            \ifnum \i=\tnk
                \node[left=of fine_\i] (tnk) {\( t^{n,k} \)};
                \node[above=-1pt of tnk] {\( \vdots \)};
                \node[below=-8pt of tnk] {\( \vdots \)};
                \draw[thick,dashed] (0,\height*\i/\numfine) -- (fine_\i);
            \else
                \draw[thick] (0,\height*\i/\numfine) -- (fine_\i);
                \ifnum \i=\inteval{\numfine-1}
                    \node[left=of fine_\i] {\( t^{n,M-1} \)};
                \else
                    \ifnum \i=\numfine
                        \node[left=of fine_\i] {\( t^{n,M} \)};
                    \else
                        \node[left=of fine_\i] {\( t^{n,\i} \)};
                    \fi
                \fi
            \fi
        }

        \node (fine_first) at (-\width/1.25,\height*\tnk/\numfine+\height/\numfine/3) {};
        \draw[thick] (0,\height*\tnk/\numfine+\height/\numfine/3) -- (fine_first);
        \node[below left=-15pt of fine_first] {\( t^{n,k+\frac{1}{3}} \)};

        \node (fine_second) at (-\width/1.25,\height*\tnk/\numfine+\height/\numfine/2) {};
        \draw[thick] (0,\height*\tnk/\numfine+\height/\numfine/2) -- (fine_second);
        \node[above left=-15pt of fine_second] {\( t^{n,k+\frac{1}{2}} \)};

        \draw[decoration={brace,raise=0pt},decorate,thick] (-\width-3,0) -- node[left=5pt] {\( \frac{\Delta t}{M} \)} (-\width-3,\height/\numfine);
    \end{tikzpicture}
    \caption{
        Visualization and notation for the time-levels used by FB-LTS.
        In the coarse region (illustrated in the right half of the diagram), \( t^{n+\nicefrac{1}{m}} = t^n + \frac{\Delta t}{m} \) for any positive integer \( m \); first stage data is calculated at time \( t^{n+\nicefrac{1}{3}} \) and second stage data is calculated at time \( t^{n+\nicefrac{1}{2}} \).
        In the fine region (illustrated in the left half of the diagram), \( t^{n,k} = t^n + k\frac{\Delta t}{M} \) and \( t^{n,k+\nicefrac{1}{m}} = t^n + \left(k+\nicefrac{1}{m} \right)\frac{\Delta t}{M} \); first stage data is calculated at times \( t^{n,k+\nicefrac{1}{3}} \) and second stage data is calculated at times \( t^{n,k+\nicefrac{1}{2}} \) for \( k = 0,\cdots, M-1 \).
    }
    \label{fig:time_levels}
\end{figure}
The FB-LTS scheme proceeds as follows.
\begin{enumerate}
    \item \textbf{Coarse Advancement:} Compute all three stages of FB-RK(3,2) on cells and edges from \( \C_P \) and \( \C_E \) to advance to time \( t^{n+1} \).
    Note that during this step, we perform calculations on some of the interface adjacent fine cells; recall that if \( r \) is the radius of the discrete tendency operators for \( \u \) and \( \h \), we defined \( \F_P^\ell \) for \( \ell = 1, \cdots, 5 \) so that \( \F_P^\ell \) contains \( \ell r \) layers of cells in \( \F_P \) neighboring \( \C_P^\text{IF-1} \), and that \( \F_E^\ell \) was the corresponding sets of edges for \( \ell = 1, \cdots, 5 \).
    The data computed in these sets is not used to advance the fine region, but is needed because these cells and edges are in the domain of dependence for the interface regions.
    \begin{enumerate}
        \item \textbf{Thickness Stage 1:} For \( i \in \F_P^5 \cup \C_P^\ifone \cup \C_P^\iftwo \),
        \begin{linenomath}
        \begin{equation}
        \begin{aligned}
            \tilde{h}_i^{n+\nicefrac{1}{3}} &= h^n + \frac{\Delta t}{3} \Psi_i\left( \u^n, \h^n \right) \\
            h_i^* &= \beta_1 \tilde{h}_i^{n+\nicefrac{1}{3}} + (1-\beta_1) h_i^n \,.
        \end{aligned} \label{eqn:if_h_s1}
        \end{equation}
        \end{linenomath}
        Note that \( \tilde{h}_i \) denotes what we refer to as \emph{uncorrected} thickness data. At the end of the scheme we will recalculate the necessary uncorrected data on the interface regions using information from the fine region advancement to obtain \emph{corrected} data..

        For \( i \in \C_P^\interior \),
        \begin{linenomath}
        \begin{equation}
        \begin{aligned}
            \bar{h}_i^{n+\nicefrac{1}{3}} &= h^n + \frac{\Delta t}{3} \Psi_i\left( \u^n, \h^n \right) \\
            h_i^* &= \beta_1 \bar{h}_i^{n+\nicefrac{1}{3}} + (1-\beta_1) h_i^n \,.
        \end{aligned} \label{eqn:coarse_h_s1}
        \end{equation}
        \end{linenomath}

        \item \textbf{Velocity Stage 1:} For \( e \in \F_E^4 \cup \C_E^\ifone \cup \C_E^\iftwo \),
        \begin{linenomath}
        \begin{equation}
        \begin{aligned}
            \tilde{u}_e^{n+\nicefrac{1}{3}} &= u^n + \frac{\Delta t}{3} \Phi_e\left( \u^n, \h^* \right) \,.
        \end{aligned} \label{eqn:if_u_s1}
        \end{equation}
        \end{linenomath}
        Note that \( \tilde{u}_e \) denotes what we refer to as \emph{uncorrected} velocity data. At the end of the scheme we will recalculate the necessary uncorrected data on the interface regions using information from the fine region advancement to obtain \emph{corrected} data.

        For \( e \in \C_E^\interior \),
        \begin{linenomath}
        \begin{equation}
        \begin{aligned}
            \tilde{u}_e^{n+\nicefrac{1}{3}} &= u^n + \frac{\Delta t}{3} \Phi_e\left( \u^n, \h^* \right) \,.
        \end{aligned} \label{eqn:coarse_u_s1}
        \end{equation}
        \end{linenomath}

        \item \textbf{Thickness Stage 2:} For \( i \in \F_P^3 \cup \C_P^\ifone \cup \C_P^\iftwo \),
        \begin{linenomath}
        \begin{equation}
        \begin{aligned}
            \tilde{h}_i^{n+\nicefrac{1}{2}} &= h_i^n + \frac{\Delta t}{2} \Psi_i\left( \tilde{\u}^{n+\nicefrac{1}{3}}, \tilde{\h}^{n+\nicefrac{1}{3}} \right) \\ 
            h_i^{**} &= \beta_2 \tilde{h}_i^{n+\nicefrac{1}{2}} + (1-\beta_2) h_i^n\,,
        \end{aligned} \label{eqn:if_h_s2}
        \end{equation}
        \end{linenomath}
        where \( \tilde{h}_i^{n+\nicefrac{1}{3}} := \bar{h}_i^{n+\nicefrac{1}{3}} \) and \( \tilde{u}_e^{n+\nicefrac{1}{3}} := \bar{u}_e^{n+\nicefrac{1}{3}} \) for \( i \in \C_P^\interior \) and \( e \in \C_E^\interior \).

        For \( i \in \C_P^\interior \),
        \begin{linenomath}
        \begin{equation}
        \begin{aligned}
            \bar{h}_i^{n+\nicefrac{1}{2}} &= h_i^n + \frac{\Delta t}{2} \Psi_i\left( \bar{\u}^{n+\nicefrac{1}{3}}, \bar{\h}^{n+\nicefrac{1}{3}} \right) \\ 
            h_i^{**} &= \beta_2 \bar{h}_i^{n+\nicefrac{1}{2}} + (1-\beta_2) h_i^n\,,
        \end{aligned} \label{eqn:coarse_h_s2}
        \end{equation}
        \end{linenomath}
        where \( \bar{h}_i^{n+\nicefrac{1}{3}} := \tilde{h}_i^{n+\nicefrac{1}{3}} \) and \( \bar{u}_e^{n+\nicefrac{1}{3}} := \tilde{u}_e^{n+\nicefrac{1}{3}} \) for \( i \in \C_P^\iftwo \) and \( e \in \C_E^\iftwo \).

        \item \textbf{Velocity Stage 2:} For \( e \in \F_E^2 \cup \C_E^\ifone \cup \C_E^\iftwo \),
        \begin{linenomath}
        \begin{equation}
        \begin{aligned}
            \tilde{u}_e^{n+\nicefrac{1}{2}} = u_e^n + \frac{\Delta t}{2} \Phi_e\left( \tilde{\u}^{n+\nicefrac{1}{3}}, \h^{**} \right)\,,
        \end{aligned} \label{eqn:if_u_s2}
        \end{equation}
        \end{linenomath}
        where \( \tilde{u}_e^{n+\nicefrac{1}{3}} := \bar{u}_e^{n+\nicefrac{1}{3}} \) for \( e \in \C_E^\interior \).

        For \( e \in \C_E^\interior \),
        \begin{linenomath}
        \begin{equation}
        \begin{aligned}
            \bar{u}_e^{n+\nicefrac{1}{2}} = u_e^n + \frac{\Delta t}{2} \Phi_e\left( \bar{\u}^{n+\nicefrac{1}{3}}, \h^{**} \right)\,,
        \end{aligned} \label{eqn:coarse_u_s2}
        \end{equation}
        \end{linenomath}
        where \( \bar{u}_e^{n+\nicefrac{1}{3}} := \tilde{u}_e^{n+\nicefrac{1}{3}} \) for \( e \in \C_E^\iftwo \).

        \item \textbf{Thickness Stage 3:} For \( i \in \F_P^1 \cup \C_P^\ifone \cup \C_P^\iftwo \),
        \begin{linenomath}
        \begin{equation}
        \begin{aligned}
            \tilde{h}_i^{n+1} &= h_i^n + \Delta t \Psi_i\left( \tilde{\u}^{n+\nicefrac{1}{2}}, \tilde{\h}^{n+\nicefrac{1}{2}} \right) \\ 
            h_i^{***} &= \beta_3 \tilde{h}_i^{n+1} + (1-2\beta_3) \tilde{h}_i^{n+\nicefrac{1}{2}} + \beta_3 h_i^n\,,
        \end{aligned} \label{eqn:if_h_s3}
        \end{equation}
        \end{linenomath}
        where \( \tilde{h}_i^{n+\nicefrac{1}{2}} := \bar{h}_i^{n+\nicefrac{1}{2}} \) and \( \tilde{u}_e^{n+\nicefrac{1}{2}} := \bar{u}_e^{n+\nicefrac{1}{2}} \) for \( i \in \C_P^\interior \) and \( e \in \C_E^\interior \).

        For \( i \in \C_P^\interior \),
        \begin{linenomath}
        \begin{equation}
        \begin{aligned}
            h_i^{n+1} &= h_i^n + \Delta t \Psi_i\left( \bar{\u}^{n+\nicefrac{1}{2}}, \bar{\h}^{n+\nicefrac{1}{2}} \right) \\ 
            h_i^{***} &= \beta_3 h_i^{n+1} + (1-2\beta_3) \bar{h}_i^{n+\nicefrac{1}{2}} + \beta_3 h_i^n\,,
        \end{aligned} \label{eqn:coarse_h_s3}
        \end{equation}
        \end{linenomath}
        where \( \bar{h}_i^{n+\nicefrac{1}{2}} := \tilde{h}_i^{n+\nicefrac{1}{2}} \) and \( \bar{u}_e^{n+\nicefrac{1}{2}} := \tilde{u}_e^{n+\nicefrac{1}{2}} \) for \( i \in \C_P^\iftwo \) and \( e \in \C_E^\iftwo \).

        \item \textbf{Velocity Stage 3:} For \( e \in \C_E^\ifone\),
        \begin{linenomath}
        \begin{equation}
        \begin{aligned}
            \tilde{u}_e^{n+1} = u_e^n + \Delta t \Phi_e\left( \tilde{\u}^{n+\nicefrac{1}{2}}, \h^{***} \right)\,,
        \end{aligned}
        \end{equation}
        \end{linenomath}
        where \( \tilde{u}_e^{n+\nicefrac{1}{2}} := \bar{u}_e^{n+\nicefrac{1}{2}} \) for \( e \in \C_E^\interior \).

        For \( e \in \C_E^\interior \),
        \begin{linenomath}
        \begin{equation}
        \begin{aligned}
            u_e^{n+1} = u_e^n + \Delta t \Phi_e\left( \bar{\u}^{n+\nicefrac{1}{2}}, \h^{***} \right)\,,
        \end{aligned}
        \end{equation}
        \end{linenomath}
        where \( \bar{u}_e^{n+\nicefrac{1}{2}} := \tilde{u}_e^{n+\nicefrac{1}{2}} \) for \( e \in \C_E^\iftwo \).
    \end{enumerate}

    \item \textbf{Interface Prediction:} Use the uncorrected data on interface one to obtain predicted values for FB-RK(3,2) data on \( \C_{P}^\ifone \) and \( \C_E^\ifone \) at times \( t^{n,k} \), \( t^{n,k+\nicefrac{1}{3}} \), and \( t^{n,k+\nicefrac{1}{2}} \), for \( k = 0, \cdots, M-1 \).
    \begin{linenomath}
    \begin{subequations}
    \begin{align}
        \begin{bmatrix}
            \h^{n,k} \\
            \u^{n,k}
        \end{bmatrix} &= \frac{k}{M}\begin{bmatrix}
            \tilde{\h}^{n+1} \\
            \tilde{\u}^{n+1}
        \end{bmatrix} + \left(1-\frac{k}{M}\right)\begin{bmatrix}
            \h^n \\
            \u^n
        \end{bmatrix} \label{subeqn:pred_old} \\
        &\null \nonumber \\
        \begin{bmatrix}
            \bar{\h}^{n,k+\nicefrac{1}{3}} \\
            \bar{\u}^{n,k+\nicefrac{1}{3}}
        \end{bmatrix} &= \frac{k}{M}\begin{bmatrix}
            \tilde{\h}^{n+1} \\
            \tilde{\u}^{n+1}
        \end{bmatrix} + \frac{1}{M}\begin{bmatrix}
            \tilde{\h}^{n+\nicefrac{1}{3}} \\
            \tilde{\u}^{n+\nicefrac{1}{3}}
        \end{bmatrix} + \left(1-\frac{k+1}{M}\right)\begin{bmatrix}
            \h^n \\
            \u^n
        \end{bmatrix} \label{subeqn:pred_s1} \\
        &\null \nonumber \\
        \begin{bmatrix}
            \bar{\h}^{n,k+\nicefrac{1}{2}} \\
            \bar{\u}^{n,k+\nicefrac{1}{2}}
        \end{bmatrix} &= \frac{k}{M}\begin{bmatrix}
            \tilde{\h}^{n+1} \\
            \tilde{\u}^{n+1}
        \end{bmatrix} + \frac{1}{M}\begin{bmatrix}
            \tilde{\h}^{n+\nicefrac{1}{2}} \\
            \tilde{\u}^{n+\nicefrac{1}{2}}
        \end{bmatrix} + \left(1-\frac{k+1}{M}\right)\begin{bmatrix}
            \h^n \\
            \u^n
        \end{bmatrix}\,. \label{subeqn:pred_s2}
    \end{align} \label{eqn:preds}
    \end{subequations}
    \end{linenomath}
    Additionally, compute the prediction for \( \h^{n,k} \) \eqref{subeqn:pred_old} one additional time for \( k = M \); this data will be needed to calculate the third stage of FB-RK(3,2) in the fine region with \( k = M-1 \).
    The coefficients for interpolating the uncorrected data on interface one are called the prediction coefficients, and the chosen prediction coefficients given above result in a second order approximation to \( \u \) and \( \h \) data at the corresponding times.
    These coefficients are derived in section \ref{sec:derivation_pred_coeffs}.
    Additionally, one can observe that when \( M = 1 \), \eqref{eqn:preds} reduces in such a way that the predictions are exactly the corresponding data already computed.
    This means that in the case where \( M = 1 \), FB-LTS is mathematically equivalent to FB-RK(3,2).

    Note that we can also use this data to compute values for \( \h^{*,k} \), \( \h^{**,k} \), and \( \h^{***,k} \) on interface one cells. For \( i \in \C_P^\ifone \),
    \begin{linenomath}
    \begin{equation}
    \begin{split}
        h_i^{*,k} &= \beta_1 \bar{h}_i^{n,k+\nicefrac{1}{3}} + (1-\beta_1) h_i^{n,k} \\
        h_i^{**,k} &= \beta_2 \bar{h}_i^{n,k+\nicefrac{1}{2}} + (1-\beta_2) h_i^{n,k} \\
        h_i^{***,k} &= \beta_3 h_i^{n,k+1} + (1-2\beta_3) \bar{h}_i^{n,k+\nicefrac{1}{2}} + \beta_3 h_i^{n,k} \,.
    \end{split} \label{eqn:hstar_preds}
    \end{equation}
    \end{linenomath}

    \item \textbf{Fine Advancement:} For \( i \in \F_P \) and \( e \in \F_E \), advance with FB-RK(3,2) with the fine time-step \( M \) times. 
    For \( k = 0, \cdots, M-1 \),
    \begin{linenomath}
    \begin{subequations}
    \label{eqn:fine_fbrk32}
    \begin{align}
    \begin{split}
        \bar{h}_i^{n,k+\nicefrac{1}{3}} &= h_i^{n,k} + \frac{\Delta t}{3M} \Psi_i\left( \u^{n,k}, \h^{n,k} \right) \\
        \bar{u}_e^{n,k+\nicefrac{1}{3}} &= u_e^{n,k} + \frac{\Delta t}{3M} \Phi_e\left( \u^{n,k}, \h^{*,k} \right) \\
        h_i^{*,k} &= \beta_1 \bar{h}_i^{n,k+\nicefrac{1}{3}} + (1-\beta_1) h_i^{n,k}
    \end{split} \label{subeqn:fine_s1} \\
    \nonumber \\
    \begin{split}
        \bar{h}_i^{n,k+\nicefrac{1}{2}} &= h_i^{n,k} + \frac{\Delta t}{2M} \Psi_i\left( \bar{\u}^{n,k+\nicefrac{1}{3}}, \bar{\h}^{n,k+\nicefrac{1}{3}} \right) \\
        \bar{u}_e^{n,k+\nicefrac{1}{2}} &= u_e^{n,k} + \frac{\Delta t}{2M} \Phi_e\left( \bar{\u}^{n,k+\nicefrac{1}{3}}, \h^{**,k} \right) \\
        h_i^{**,k} &= \beta_2 \bar{h}_i^{n,k+\nicefrac{1}{2}} + (1-\beta_2) h_i^{n,k}
    \end{split} \label{subeqn:fine_s2} \\
    \nonumber \\
    \begin{split}
        h_i^{n,k+1} &= h_i^{n,k} + \frac{\Delta t}{M} \Psi_i\left( \bar{\u}^{n,k+\nicefrac{1}{2}}, \bar{\h}^{n,k+\nicefrac{1}{2}} \right) \\
        u_e^{n,k+1} &= u_e^{n,k} + \frac{\Delta t}{M} \Phi_e\left( \bar{\u}^{n,k+\nicefrac{1}{2}}, \h^{***,k} \right) \\
        h_i^{***,k} &= \beta_3 h_i^{n,k+1} + (1-2\beta_3) \bar{h}_i^{n,k+\nicefrac{1}{2}} + \beta_3 h_i^{n,k} \,.
    \end{split} \label{subeqn:fine_s3}
    \end{align}
    \end{subequations}
    \end{linenomath}
    Note that in the FB average in \eqref{subeqn:fine_s3}, we already have the data for \( h_i^{n,k+1} \) on \( \C^\ifone \) cells for \( k = M-1 \) because of the extra computation of \eqref{subeqn:pred_old} for \( k = M \).

    After looping over \( k \), set \( h_i^{n+1} := h_i^{n,M} \) and \( u_e^{n+1} := u_e^{n,M} \) for \( i \in \F_P \) and \( e \in \F_E \).

    \item \textbf{Interface Correction:} Finally, compute the corrected data at time \( t^{n+1} \) on interface one and two. 
    For \( i \in \C_P^\ifone \cup \C_P^\iftwo \) and \( e \in \C_E^\ifone \cup \C_E^\iftwo \),
    \begin{linenomath}
    \begin{equation}
    \begin{aligned}
        h_i^{n+1} &= h_i^{n} + \frac{\Delta t}{M} \sum_{k=0}^{M-1} \Psi_i\left( \bar{\u}^{n,k+\nicefrac{1}{2}}, \bar{\h}^{n,k+\nicefrac{1}{2}} \right) \\
        u_e^{n+1} &= u_e^{n} + \frac{\Delta t}{M} \sum_{k=0}^{M-1} \Phi_e\left( \bar{\u}^{n,k+\nicefrac{1}{2}}, \h^{***,k} \right)\,,
    \end{aligned} \label{eqn:if_correction}
    \end{equation}
    \end{linenomath}
    where  \( u_e^{n,k} := u_e^n \), \( h_i^{n,k} := h_i^n \), \( h_i^{*,k} := h_i^* \), \( \bar{u}_e^{n,k+\nicefrac{1}{3}} := \tilde{u}_e^{n+\nicefrac{1}{3}} \), \( \bar{h}_i^{n,k+\nicefrac{1}{3}} := \tilde{h}_i^{n+\nicefrac{1}{3}} \), \( h_i^{**,k} := h_i^{**} \), \( \bar{u}_e^{n,k+\nicefrac{1}{2}} := \tilde{u}_e^{n+\nicefrac{1}{2}} \), \( \bar{h}_i^{n,k+\nicefrac{1}{2}} := \tilde{h}_i^{n+\nicefrac{1}{2}} \), and \( h_i^{***,k} := h_i^{***} \) for \( i \in \C_P^\iftwo \) and \( e \in \C_E^\iftwo \).
    Similarly, \( u_e^{n,k} := u_e^n \), \( h_i^{n,k} := h_i^n \), \( h_i^{*,k} := h_i^* \), \( \bar{u}_e^{n,k+\nicefrac{1}{3}} := \bar{u}_e^{n+\nicefrac{1}{3}} \), \( \bar{h}_i^{n,k+\nicefrac{1}{3}} := \bar{h}_i^{n+\nicefrac{1}{3}} \), \( h_i^{**,k} := h_i^{**} \), \( \bar{u}_e^{n,k+\nicefrac{1}{2}} := \bar{u}_e^{n+\nicefrac{1}{2}} \), \( \bar{h}_i^{n,k+\nicefrac{1}{2}} := \bar{h}_i^{n+\nicefrac{1}{2}} \), and \( h_i^{***,k} := h_i^{***} \) for \( i \in \C_P^\interior \) and \( e \in \C_E^\interior \).

    Note that the individual terms in the sums in \eqref{eqn:if_correction} can be calculated and accumulated during the fine advancement step.
\end{enumerate}

This ends the description of the FB-LTS scheme.


\subsection{Temporal Convergence}
\label{subsec:convergence}

Here we describe the results of a numerical experiment in which FB-LTS is \( \Ord{2} \) everywhere, including on interface cells and edges (Figure \ref{fig:convergence}).
We calculate the error of a FB-LTS solution against that of the classical, four-stage, fourth order Runge-Kutta method (RK4) using a small time-step on a simple model problem.

Consider a non-rotating aquaplanet (i.e. a spherical mesh with no land cells) with constant resting thickness, where the fluid velocity is initialized to zero, and the layer thickness is initialized to a Gaussian bump.
This produces a simple external gravity wave, which is modeled by the  system
\begin{linenomath}
\begin{equation}
    \begin{cases}
        \pd{\u}{t} = - g\nabla h \\
        \pd{h}{t} + \nabla \cdot (h \mathbf{u}) = 0 \,.
    \end{cases}
    \label{eqn:grav_wave_model}
\end{equation}
\end{linenomath}
In these equations
\( \mathbf{u} \) is the fluid velocity,
\( t \) is the time coordinate, 
\( g \) is the gravitational constant, 
and \( h \) is the ocean thickness.

The root-mean-square (RMS) error is defined as
\begin{linenomath}
\begin{equation}
    E_{\text{RMS}} = \sqrt{\frac{\sum_{i = 1}^N (s_i - m_i)^2}{N}} \,,
    \label{eqn:rmse}
\end{equation}
\end{linenomath}
where \( \{ s_i \}_{i=1}^N \) is the discrete reference solution, \( \{m_i\}_{i=1}^N \) is the discrete model solution, and \( N \) is the number of discretization points.

\begin{figure}
    \centering
    \includegraphics[width=0.8\textwidth]{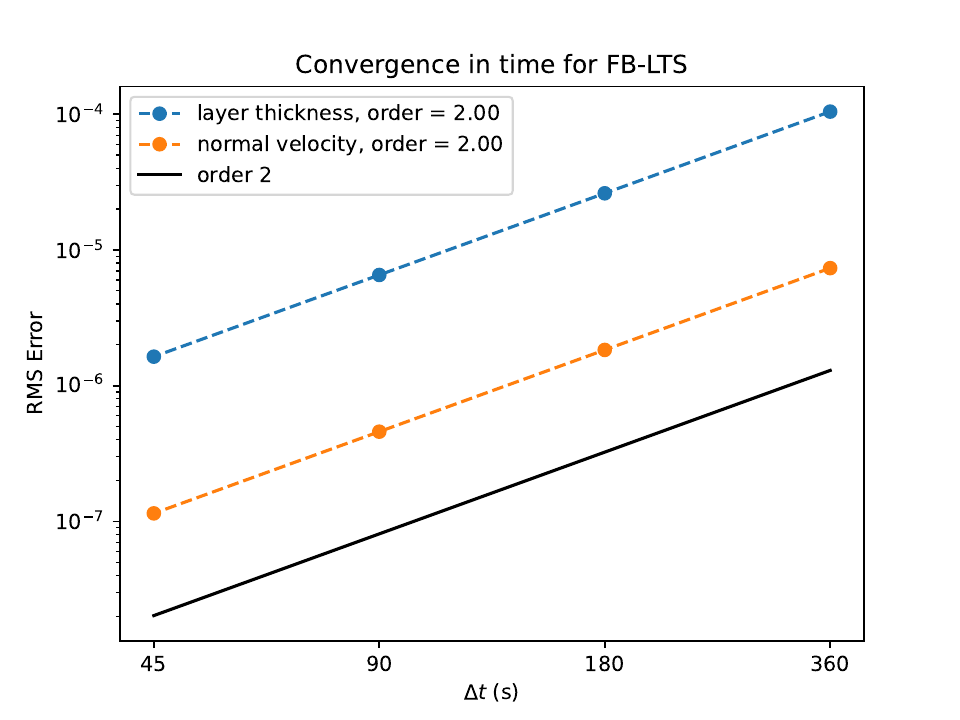}
    \caption{
        Temporal convergence for FB-LTS on a test case consisting of an external gravity wave on a non-rotating aquaplanet.
        The time-steps on the horizontal axis are the time-steps used in the coarse region, in each case we have \( M = 4 \).
        Errors are computed against a reference solution generated by RK4 using a time-step of 10 s.
    }
    \label{fig:convergence}
\end{figure}


\subsection{Conservation of Mass and Absolute Vorticity}
\label{subsec:conservation}

An important property of FB-LTS is that it provides exact conservation of the spatially discrete representation of mass (Theorem \ref{thm:conservation_of_mass}) and the spatially discrete representation of absolute vorticity (Theorem \ref{thm:conservation_of_abs_vorticity}) in the cases where there are no boundary conditions (e.g. the case of an aquaplanet), or no-flow boundary conditions (i.e. the normal velocity is zero at the boundary).
In the context of the shallow water equations, the mass variable is the ocean layer thickness \( h \), and absolute vorticity is given by \( \eta = \mathbf{k} \cdot \nabla \times \u + f \), depending on the fluid velocity \( \u \) and the Coriolis parameter \( f \).
To clarify what is meant by conservation here, we mean that the values of globally integrated thickness and globally integrated absolute vorticity are constant in time.
We are concerned with the discrete counterparts to the conservation equations for these quantities in the continuous case.
The mass equation is given by
\begin{linenomath}
\begin{equation}
    \pd{h}{t} + \nabla \cdot (h\u) = 0 \,,
    \label{eqn:mass_conservation}
\end{equation}
\end{linenomath}
and the absolute vorticity equation (obtained by taking the curl of the momentum equation) is given by
\begin{linenomath}
\begin{equation}
    \pd{\eta}{t} + \nabla \cdot \left( \eta \u \right) = 0 \,.
    \label{eqn:abs_vorticity_conservation}
\end{equation}
\end{linenomath}
Often, \eqref{eqn:abs_vorticity_conservation} is written in terms of potential vorticity (PV), where PV is given by \( q = \frac{\eta}{h} \).
The absolute vorticity equation can also be written as a thickness-weighted potential vorticity equation,
\begin{linenomath}
\begin{equation}
    \pd{(qh)}{t} + \nabla \cdot \left(  q(h\u) \right) = 0 \,.
    \label{eqn:thick_pv_conservation}
\end{equation}
\end{linenomath}
In particular, \citet{ringler2010} formulates the discrete counterpart to \eqref{eqn:abs_vorticity_conservation} in terms of a spatially discrete representation of PV.

\begin{theorem}\label{thm:conservation_of_mass}
    FB-LTS exactly preserves the discrete representation of mass assuming either no boundary conditions or no-flow boundary conditions.
\end{theorem}

\begin{proof}
    The continuous equation giving conservation of mass for the SWEs is \eqref{eqn:mass_conservation},
    which states that the evolution of the thickness \( h \) depends only on the divergence of the thickness flux.
    Within a TRiSK framework, this equation is spatially discretized as
    \begin{linenomath}
    \begin{align}
        \pd{h_i}{t} &= \Psi_i\left( \u, \h \right) \nonumber \\
        &= \frac{-1}{A_i} \sum_{e \in \mathcal{E}_i} n_{e,i} \ell_e F_e\left(u_e, \h \right) \,,
        \label{eqn:trisk_thickness}
    \end{align}
    \end{linenomath}
    where \( \u = (u_e)_{e \in \Omega_E} \) and \( \h = (h_i)_{i \in \Omega_P} \),
    \( A_i \) is the area of cell \( i \),
    \( \mathcal{E}_i \) is the set of edges of cell \( i \),
    \( \ell_e \) is the length of edge \( e \),
    \( n_{e,i} \) is either \( 1 \) or \( -1 \), chosen so that \( n_{e,i} \mathbf{n}_e \) (see Figure \ref{fig:trisk}) is the unit outward normal vector to cell \( i \) at edge \( e \),
    and \( F_e =  F_e\left(u_e, \h \right) \) is the signed value of the thickness flux at edge \( e \) in the direction of the unit vector \( \mathbf{n}_e \).
    Note that because of the choice of \( n_{e,i} \) such that \( n_{e,i}\mathbf{n}_e \) is the outward normal to cell \( i \) at edge \( e \), the quantity \( n_{e,i} F_e \) is the thickness flux leaving cell \( i \).
    To show that the total mass in the system at time \( t^n \) is equal to the total mass in the system at time \( t^{n+1} \), we will show that any given edge is neither a source nor a sink for mass as it is transported across cells.
    Specifically, we will show this for an edge that is shared by a fine cell and an interface one cell as this is the most delicate case; the cases of other types of edges can be shown similarly, though more simply, so we omit them for brevity.

    Let \(i_1 \in \F_P \) and \( i_2 \in \C_P^\ifone \) be indices for cells \( P_{i_1}\) and \( P_{i_2} \) that have a shared edge \( e_0 \in \F_E \).
    Assume that the thickness-fluxes at all other edges of \( P_{i_1} \) and \( P_{i_2} \) are zero.
    We do this with the goal of showing that the mass flux as computed on one side of the edge \( e_0 \) is the same as computed on the other side of the edge \( e_0 \).
    From this, we can conclude that the edge in question neither adds or subtracts mass from the system, giving a local conservation of mass, and from there, global conservation follows.
    Since we are assuming that \( F_e\left(u_e, \h \right) = 0 \) for all \( e \neq e_0 \), \eqref{eqn:trisk_thickness} simplifies to
    \begin{linenomath}
    \begin{equation}
        \pd{h_i}{t} = \frac{-1}{A_i} n_{e_0,i} \ell_{e_0} F_{e_0}\left(u_{e_0}, \left(h_{i_1}, h_{i_2}\right) \right) \,,
        \label{eqn:trisk_thickness_simplify}
    \end{equation}
    \end{linenomath}
    for \( i = i_1, i_2 \).
    The goal is to show that
    \begin{linenomath}
    \begin{equation}
        \int_{P_{i_1} \cup P_{i_2}} h_{i_*}^{n+1} \dif A = \int_{P_{i_1} \cup P_{i_2}} h_{i_*}^n \dif A \,,
        \label{eqn:conservation}
    \end{equation}
    \end{linenomath}
    where  \( h_{i_*} \) refers to the thicknesses \( h_{i_1} \) and \( h_{i_2} \) in the respective cells.
    From here, the central idea is to look at the discretization of the thickness on cell \( i_1 \) and \( i_2 \) at time \( t^{n+1} \).
    For cell \( i_1 \), which is in the fine region, the thickness at time \( t^{n+1} \) is given by \eqref{subeqn:fine_s3} with \( k = M-1 \) in the \textit{Fine Advancement} step of FB-LTS.
    For cell \( i_2 \), which is in the interface one region, the thickness at time \( t^{n+1} \) is given by \eqref{eqn:if_correction} in the \textit{Interface Correction} step.
    After writing out both \( h_{i_1}^{n+1} \) and \( h_{i_2}^{n+1} \), we will note that the total flux over their shared edge is equal, as computed from both cells.
    Then, we can integrate both thicknesses over their respective cells to show that the total mass in both cells at time \( t^{n+1} \) is the same as the total mass at time \( t^n \).
    On \( P_{i_1} \), we have that
    \begin{linenomath}
    \begin{align}
        h_{i_1}^{n+1} &= h_{i_1}^{n,M} \nonumber \\
        &= h_{i_1}^{n,M-1} + \frac{\Delta t}{M} \Psi_{i_1}\left( \bar{\u}^{n,(M-1)+\nicefrac{1}{2}}, \bar{\h}^{n,(M-1)+\nicefrac{1}{2}} \right) \nonumber \\
        &= h_{i_1}^{n,M-2} + \frac{\Delta t}{M}\left[ \Psi_{i_1}\left( \bar{\u}^{n,(M-2)+\nicefrac{1}{2}}, \bar{\h}^{n,(M-2)+\nicefrac{1}{2}} \right) + \Psi_{i_1}\left( \bar{\u}^{n,(M-1)+\nicefrac{1}{2}}, \bar{\h}^{n,(M-1)+\nicefrac{1}{2}} \right) \right] \nonumber \\
        &\vdots \nonumber \\
        &= h_{i_1}^{n,0} + \frac{\Delta t}{M} \sum_{k=0}^{M-1} \Psi_{i_1}\left( \bar{\u}^{n,k+\nicefrac{1}{2}}, \bar{\h}^{n,k+\nicefrac{1}{2}} \right) \nonumber \\
        &= h_{i_1}^n + \frac{\Delta t}{M} \sum_{k=0}^{M-1}\left[ \frac{-1}{A_{i_1}} n_{e_0,{i_1}} \ell_{e_0} F_{e_0}\left(\bar{u}_{e_0}^{n,k+\nicefrac{1}{2}}, \left(\bar{h}_{i_1}^{n,k+\nicefrac{1}{2}}, \bar{h}^{n,k+\nicefrac{1}{2}}_{i_2}\right) \right) \right] \,. \label{eqn:h_P1}
    \end{align}
    \end{linenomath}
    Note that the values of the thickness fluxes \( F_{e_0}\left(\bar{u}_{e_0}^{n,k+\nicefrac{1}{2}}, \left(\bar{h}_{i_1}^{n,k+\nicefrac{1}{2}}, \bar{h}^{n,k+\nicefrac{1}{2}}_{i_2}\right) \right) \) depend only on the thickness on cell \( i_1 \) in the fine region, the thickness on cell \( i_2 \) in the interface one region, and the normal velocity at edge \( e_0 \), all computed as second stage data times \( t^{n,k+\nicefrac{1}{2}} \) for \( k = 0, \cdots, M-1 \).
    To simplify notation, from here on we write \( F_{e_0}^{n,k+\nicefrac{1}{2}} := F_{e_0}\left(\bar{u}_{e_0}^{n,k+\nicefrac{1}{2}}, \left(\bar{h}_{i_1}^{n,k+\nicefrac{1}{2}}, \bar{h}^{n,k+\nicefrac{1}{2}}_{i_2}\right) \right) \).
    The central idea for the proof is that both the fine cell and the interface one cell are using the same flux for each \( k = 0,\cdots, M-1 \).
    Now, looking at \( P_{i_2} \), after completing the correction step \eqref{eqn:if_correction}, we have
    \begin{linenomath}
    \begin{align}
        h_{i_2}^{n+1} &= h_{i_2}^n + \frac{\Delta t}{M} \sum_{k=0}^{M-1} \Psi_{i_2}\left( \bar{\u}^{n,k+\nicefrac{1}{2}}, \bar{\h}^{n,k+\nicefrac{1}{2}} \right) \nonumber \\
        &=  h_{i_2}^n + \frac{\Delta t}{M} \sum_{k=0}^{M-1} \left[ \frac{-1}{A_{i_2}} n_{e_0,{i_2}} \ell_{e_0} F_{e_0}^{n,k+\nicefrac{1}{2}} \right] \,. \label{eqn:h_P2}
    \end{align}
    \end{linenomath}
    Next, we integrate \eqref{eqn:h_P1} over the cell \( P_{i_1} \) to get the total mass in the cell at time \( t^{n+1} \),
    \begin{linenomath}
    \begin{align}
        \int_{P_{i_1}} h_{i_1}^{n+1} \dif A &= A_{i_1} h_{i_1}^n + A_{i_1} \left( \frac{\Delta t}{M} \sum_{k=0}^{M-1}\left[ \frac{-1}{A_{i_1}} n_{e_0,{i_1}} \ell_{e_0} F_{e_0}^{n,k+\nicefrac{1}{2}} \right] \right) \nonumber \\
        &= A_{i_1} h_{i_1}^n  - \frac{\Delta t}{M} n_{e_0,{i_1}} \ell_{e_0} \sum_{k=0}^{M-1}\left[ F_{e_0}^{n,k+\nicefrac{1}{2}} \right] \nonumber \\
        &= \int_{P_{i_1}} h_{i_1}^n \dif A  - \frac{\Delta t}{M} n_{e_0,{i_1}} \ell_{e_0} \sum_{k=0}^{M-1}\left[ F_{e_0}^{n,k+\nicefrac{1}{2}} \right] \,. \label{eqn:int_h_P1}
    \end{align}
    \end{linenomath}
    Through a similar calculation on \( P_{i_2} \) starting with \eqref{eqn:h_P2}, we get
    \begin{linenomath}
    \begin{equation}
        \int_{P_{i_2}} h_{i_2}^{n+1} \dif A = \int_{P_{i_2}} h_{i_2}^n \dif A  - \frac{\Delta t}{M} n_{e_0,{i_2}} \ell_{e_0} \sum_{k=0}^{M-1}\left[ F_{e_0}^{n,k+\nicefrac{1}{2}} \right] \,. \label{eqn:int_h_P2}
    \end{equation}
    \end{linenomath}
    Finally, add \eqref{eqn:int_h_P1} and \eqref{eqn:int_h_P2},
    \begin{linenomath}
    \begin{align}
        \int_{P_{i_1} \cup P_{i_2}} h_{i_*}^{n+1} \dif A &= \int_{P_{i_1} \cup P_{i_2}} h_{i_*}^n \dif A - \frac{\Delta t}{M} \left( n_{e_0,{i_1}} + n_{e_0,{i_2}} \right) \ell_{e_0} \sum_{k=0}^{M-1}\left[ F_{e_0}^{n,k+\nicefrac{1}{2}} \right] \nonumber \\
        &= \int_{P_{i_1} \cup P_{i_2}} h_{i_*}^n \dif A \,, \label{eqn:wwwww}
    \end{align}
    \end{linenomath}
    with the final step being achieved using the fact that \( n_{e_0,{i_1}} = -n_{e_0,{i_2}} \).
    This shows that FB-LTS exactly conserves mass under the TRiSK spatial discretization.
\end{proof}

\begin{theorem}\label{thm:conservation_of_abs_vorticity}
    FB-LTS exactly preserves the discrete representation of absolute vorticity assuming either no boundary conditions or no-flow boundary conditions.
\end{theorem}

\begin{proof}
    The thickness-weighted PV equation \eqref{eqn:thick_pv_conservation} shows that the evolution of the thickness-weighted PV (i.e., absolute vorticity) depends only on the divergence of the absolute vorticity flux. 
    Under TRiSK, this is discretized in space as
    \begin{linenomath}
    \begin{align}
        \pd{\eta_v}{t} &= \Theta_v\left( \u, \h \right) \nonumber \\
        &= \frac{-1}{A_v} \sum_{e \in \mathcal{E}_v} -t_{e,v} d_e \hat{q}_e F_e^\perp \,,
        \label{eqn:trisk_vorticity}
    \end{align}
    \end{linenomath}
    where \( A_v \) is the area of dual cell \( v \),
    \( \mathcal{E}_v \) is the set of edges of dual cell \( v \),
    \( t_{v,i} \) is an indicator function, either \( -1 \) or \( 1 \)  depending on whether the fixed unit vector \( \mathbf{t}_e = \mathbf{k} \times \mathbf{n}_e \) (Figure \ref{fig:trisk}) is outward or inward to dual cell \( v \) (note that the signs are opposite relative to \( n_{i,e} \) described above; this convention is used in \citet{ringler2010} and we keep it here for consistency),
    \( d_e \) is the length of edge \( e \),
    \( \hat{q}_e \) is the value of the PV interpolated to edge \( e \),
    and \( F_e^\perp = \mathcal{M}_e\left( \mathbf{F}(\u, \h) \right) \).
    Here, \( \mathcal{M}_e \) is the flux mapping operator defined in \citet{ringler2010} that defines a mapping between the primal flux field \( \mathbf{F} = \left(F_e\right)_{e \in \Omega_E} \) in the direction normal to primal cell edges and the dual flux field \( \mathbf{F}^\perp = \left(F_e^\perp\right)_{e \in \Omega_E} \) in the direction tangent to primal cell edges, and therefore normal to dual cell edges.
    This mapping \( \mathcal{M} \) from \citet{ringler2010} allows for the vorticity equation \eqref{eqn:abs_vorticity_conservation} to be defined on the dual mesh and be consistent with the underlying momentum and thickness equations defined on the primal cells; because of this mapping, \citet{ringler2010} asserts that a prognostically obtained (i.e., computed from equations derived from PDE) vorticity is equal to a diagnostically obtained (i.e., computed from prognostically obtained quantities) vorticity, up to machine precision.
    Finally, recall that \( F_e \) and \( F_e^\perp \) are thickness fluxes of the form \( h\u \), so a quantity of the form \( q F_e = (qh)\u = \eta\u \) can be understood as a flux of the absolute vorticity.
    
    Therefore, even though in practice we obtain values for the absolute vorticity diagnostically, never computing \eqref{eqn:trisk_vorticity} directly, our goal is to show that when \eqref{eqn:trisk_vorticity} is solved prognostically, absolute vorticity is exactly conserved in the sense that the total  absolute vorticity at time \( t^n \) is equal to the total absolute vorticity at time \( t^{n+1} \).
    This can be done very similarly to the conservation of mass argument above, except now we consider dual cells, looking at fluxes normal to dual cell edges (Figure \ref{fig:dual_cells_pv}); effectively, because this is a conservation equation where the temporal evolution of \( \eta \) depends only on a divergence of a flux, it is equivalent to the thickness case shown above.
    As before, we will show this for an edge that is shared by a fine dual cell and an interface one dual cell, though calculations for other types of edges are similar.

    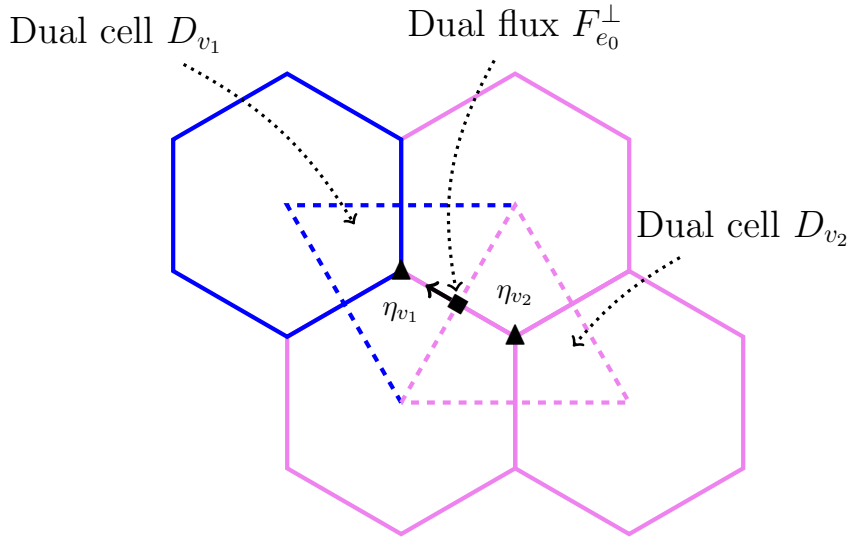
\begin{figure}
        \centering

        \def\dotscale{1.3}
        \def\mydot{\scalebox{\dotscale}{\( \mathcolor{black}{\bullet} \)}}
        \def\mysquare{\scalebox{\dotscale}{{\tiny \( \mathcolor{black}{\blacksquare} \)}}}
        \def\mytriangle{\scalebox{\dotscale}{\raisebox{1px}{\small\( \mathcolor{black}{\blacktriangle} \)}}}
    
        \tikzset{
            cell/.style={regular polygon,
            regular polygon sides=6,
            minimum size=3.5cm,
            outer sep=0, ultra thick,
            shape border rotate=30,
            text width=1.25cm,
            align=center},
        }
    
        \tikzstyle{tn}=[very thick, rounded corners=0.25cm, font=\Large, anchor=west, text width=4cm, align=center]
        \tikzstyle{ar}=[very thick, dotted, ->]
        
        \begin{tikzpicture}
            \node[cell, draw=int1outline] (int1_1) at (0,0) {};
            \node[cell, draw=int1outline, anchor=corner 1] (int1_2) at (int1_1.corner 3) {};
            \node[cell, draw=int1outline, anchor=corner 1] (int1_3) at (int1_1.corner 5) {};
            \node[cell, draw=fineoutline, anchor=corner 6] (fine) at (int1_1.corner 2) {};
            
            \draw[dashed, ultra thick, draw=fineoutline] (int1_2.center) --  (fine.center) -- (int1_1.center);
            \draw[dashed, ultra thick, draw=int1outline] (int1_2.center) -- (int1_1.center) -- (int1_3.center) -- (int1_2.center);
    
            \node (dual_1_center) at (fine.corner 5) {\mytriangle};
            \node[below=0pt of dual_1_center] {\( \eta_{v_1} \)};
    
            \node (dual_2_center) at (int1_1.corner 4) {\mytriangle};
            \node[above=0pt of dual_2_center] {\( \eta_{v_2} \)};
    
            \node[tn, above left=5pt of fine.corner 1] (dual_1_label) {Dual cell \( D_{v_1} \)};
            \coordinate[above left=25pt of fine.corner 5] (dual_1_arrow);
            \draw[ar] (dual_1_label) edge[bend left=15] (dual_1_arrow);
    
            \node[tn, above=30pt of int1_3.corner 6] (dual_2_label) {Dual cell \( D_{v_2} \)};
            \coordinate[above left=30pt of int1_3.center] (dual_2_arrow);
            \draw[ar] (dual_2_label) edge[bend right=10] (dual_2_arrow);
    
            \node[rotate=60, inner sep=0pt] (flux) at (int1_1.side 3) {\( \mysquare \)};
            \draw[->, ultra thick] (flux) -- ( $(flux)+(150:0.5)$ );
    
            \node[tn, above=5pt of int1_1.corner 1] (flux_label) {Dual flux \( F_{e_0}^\perp \)};
            \draw[ar] (flux_label) edge[bend right=20] (flux);
            
        \end{tikzpicture}
        \caption{
            TRiSK grid showing dual cells at the boundary between fine and interface one regions.
            Blue denotes the fine region and pink denotes the interface one region as in Figure \ref{fig:lts_regions}.
        }
        \label{fig:dual_cells_pv}
    \end{figure}

    Let \( \F_D \) and \( \C_D^\ifone \) be the sets of indices for dual cells, then take \( v_1 \in \F_D \) and \( v_2 \in \C_D^\ifone \) so that dual cells \( D_{v_1} \) and \( D_{v_2} \) share a dual cell edge that is orthogonal to a primal cells edge \( e_0 \in \C_E^\ifone \).
    Assume that the thickness-fluxes at all other edges of \( D_{i_1} \) and \( D_{v_2} \) are zero.
    Since we are assuming that \( F_e^\perp = 0 \) for all \( e \neq e_0 \), \eqref{eqn:trisk_vorticity} simplifies to
    \begin{linenomath}
    \begin{equation}
        \pd{\eta_v}{t} = \frac{1}{A_v} t_{e_0,v} d_{e_0} \hat{q}_{e_0} F_{e_0}^\perp \,.
        \label{eqn:trisk_vorticity_simplify}
    \end{equation}
    \end{linenomath}
    See Figure \ref{fig:dual_cells_pv} for an illustration of this situation and note how \eqref{eqn:trisk_vorticity_simplify} is similar to \eqref{eqn:trisk_thickness_simplify}.
    From here, the result follows from a calculation very similar to the above for the conservation of mass (we omit the details because the calculation is so similar to the above); the value of \( \eta_{v_1}^{n+1} \) is given by summing the fluxes in and out of \( D_{v_1} \) over \( M \) fine time-steps, and the value of \( \eta_{v_2}^{n+1} \) is given by performing the interface correction step similar to \eqref{eqn:if_correction}.
    Then, we use the fact that \( t_{e_0,v_1} = -t_{e_0,v_2} \) to conclude that the total flux entering \( D_{v_1} \) is equal to the total flux exiting \( D_{v_2} \).
    This shows that a given edge is neither a source or a sink for absolute vorticity and therefore that the absolute vorticity is conserved globally.
\end{proof}

\begin{remark}\label{rmk:volume_conservation_pv}
    Both \citet{ringler2010} and \citet{hoang2019} state that TRiSK and LTS3, respectively, conserve globally integrated PV.
    This requires the clarification that in this case, conservation means that the value of the \textit{volume} integral of PV is constant in time, as opposed to the \textit{area} integral as in the case for conservation of thickness and absolute vorticity.
    FB-LTS also conserves PV in this sense; this follows from Theorem \ref{thm:conservation_of_abs_vorticity} and the fact that \( q_{v_*} \)  is independent of the vertical position.
    That is, 
    \begin{linenomath}
    \begin{align}
        \int q_{v_*} \dif V &= \int h_{v_*}q_{v_*} \dif A \nonumber \\
        &= \int \eta_{v_*} \dif A \,. \label{eqn:vol_integral_pv}
    \end{align}
    \end{linenomath}
    Note that since PV is defined on dual cells, this requires that the thickness be defined on dual cells as well.
    TRiSK defines an auxiliary thickness equation defined on dual cells in such a way that the dual cell thickness and the primal cell thickness are equivalent with each other in the sense that a prognostically obtained dual cell thickness is equal to the dual cell thickness obtained diagnostically by an interpolation of the primal cell thickness up to machine precision \citep{ringler2010}.
\end{remark}

\begin{remark}\label{rmk:lagrangian_pv}
    We have shown that FB-LTS conserves discrete representations of both mass and absolute vorticity.
    In the continuous case, the conservation equations for these quantities \eqref{eqn:mass_conservation} and \eqref{eqn:abs_vorticity_conservation} can be combined to obtain an important statement about the evolution of PV, that PV is constant along Lagrangian trajectories, i.e.,  
    \begin{linenomath}
    \begin{equation}
        \frac{\text{D}q}{\text{D}t} = \pd{q}{t} + \u \cdot \nabla q = 0 \,,
        \label{eqn:continuous_pv_conservation}
    \end{equation}
    \end{linenomath}
    where \( \frac{\text{D}}{\text{D}t} \) is the material derivative.
    In the discrete case, using a finite volume method like TRiSK, there is no direct analogue to \eqref{eqn:continuous_pv_conservation}, so we do not claim that we conserve PV in this way.
    However, having exact conservation of mass and absolute vorticity, the statements from which \eqref{eqn:continuous_pv_conservation} can be derived in the continuous case, suggests a type of discrete analogue for this Lagrangian behavior of PV. 
    Said another way, while there is no formal sense in which \( \frac{\text{D}q}{\text{D}t} = 0 \) in the discrete case, we have satisfied the necessary pre-requisites, going as far as we can for the discrete case.
\end{remark}


\end{document}

\section{Implementation in MPAS-Ocean}
\label{sec:implementation}

Here, we describe some particulars relating to the implementation of the FB-LTS and LTS3 codes in MPAS-Ocean and compare to the now outdated implementation of LTS3 used in \citet{lilly2023_sandy} and  \citet{capodaglio2022}.


\subsection{MPI Domain Decomposition and Parallelization}
\label{subsec:parallalization}

Balancing load and achieving efficient parallelization are non-trivial tasks when it comes to local time-stepping due to the asynchronous way in which the solution is advanced, causing an inherent load imbalance that needs to be properly addressed with ad-hoc implementation strategies.
In \citet{capodaglio2022} and \citet{lilly2023_sandy}, the authors obtained load balancing and parallel scalability by assigning to each MPI rank a well balanced number of cells from the coarse, interface, and fine LTS regions.
This means that these three regions were treated as separate sets and each was partitioned across MPI ranks in a well balanced way.
This process required the interface one and interface two regions to be augmented with additional cells in order for the MPI partition to provide at least 100 interface cells per rank, as a rule of thumb.
Since the correction terms for the interface are computed during the sub-stepping procedure to minimize storage requirements, the proper number of additional interface cells needed to be tuned for best performance, as shown in \citet{capodaglio2022} and \citet{lilly2023_sandy}.

The reason for treating the interface cells as a separate entity when it comes to load balancing was motivated by the need to address another crucial feature to guarantee computational performance, which is the ability to compute the right hand side terms only on the specific regions in which the solution will be advanced.
A limitation of the underlying MPAS-Ocean framework is that the tendency routines only calculate tendency terms globally, with no way to only perform calculations on a given subset of cells.
That is, if a given MPI rank owned both fine and coarse cells, it would not be possible to calculate tendencies only on fine cells, or only on coarse cells.
In our previous works, we overcame this limitation by assigning each MPI rank multiple memory blocks, and looping over each memory block in serial.
Each memory block would only contain one type of cell (fine, interface, or coarse), so when the tendency terms were called, the tendencies would only be calculated on those cells.
Then, each MPI rank was assigned three memory blocks, each containing fine, interface, and coarse cells respectively.
This allowed us to implement the LTS algorithms without changes to the underlying MPAS-Ocean framework.

Although convenient and effective, this procedure introduced a major issue in terms of parallel communication: Each memory block had its own halo cells on which parallel communication with other MPI ranks occurred, hence the communication cost of running with \( n \) tasks was actually that of  \( 3n \) tasks.
As outlined in Table 3 from \citet{lilly2023_sandy}, the gains of using local time-stepping over a global time-stepping method were significantly reduced due to these communication issues becoming predominant when the number of cells per MPI rank became small, with 3000 cells per rank being observed as an approximate minimum for acceptable performance.

For the current work, the LTS3 and FB-LTS algorithms have been implemented in MPAS-Ocean within the main branch of the Energy Exascale Earth System Model (E3SM) \citep{golaz2022}.
This new implementation abandons the use of multiple memory blocks per MPI rank, which leads to two primary benefits. First, the amount of parallel communication does not grow by a factor of 3, hence the benefits of using local time-stepping still persist when the number of cells per MPI rank becomes small.
Second, we no longer assign one block to each LTS region, hence we no longer need to tune the number of additional cells in the interface layers, and only use the minimum number of cells in the interface layers required by the numerical algorithms. 
Concerning the load balancing procedure, the mesh cells are divided in two sets, one containing the cells in the fine LTS region, and the other containing the coarse and interface LTS regions.
Then, these two sets are each partitioned in a well balanced way according to the number of MPI ranks. Finally, each MPI rank receives its share of cells from each of the two sets. The interface cells are not treated as a separate entity anymore on the load balancing procedure  and are considered coarse cells, since the coarse time-step is used to advance the solution on them.

Discontinuing the use of multiple memory blocks per MPI rank required us to find a computationally efficient way to compute tendency terms only on the specific LTS region on which the solution is to be advanced.
This has been achieved in part via an operator splitting approach, as described in the next section.


\subsection{Operator Splitting}
\label{subsec:splitting}

It is well known that the ocean admits dynamics on a vast range of time-scales, with the most rapid motions being up to two orders of magnitude faster than the slowest.
To account for this, climate-scale ocean models employ a barotropic/baroclinic splitting in which the fast 2D motions modeled by the vertically integrated barotropic subsystem are solved separately from the slow 3D motions in the baroclinic subsystem \citep{higdon2005}.
This fast barotropic subsystem is essentially the SWEs \eqref{eqn:nonlinear_swe}; however, even within this fast subsystem, there are still a range of time-scales at play that are often solved monolithically.
The most rapid motions in a shallow water system are due to external gravity waves, which propagate at a rate of \( \sqrt{gH} \), where \( H \) is the fluid resting depth.
This means that, assuming an average ocean depth of 4000 m, gravity waves propagate at a rate of approximately 200 m \( \text{s}^{-1} \).
This is two orders of magnitude faster than Eulerian current velocities which are at most a few meters per second.

With the goal of exploiting this range of time-scales for computational performance, we introduce a certain operator splitting for the SWEs.
We identify the external gravity-wave subsystem as \textit{fast} and the remaining system, including the nonlinear momentum advection and forcing, as \textit{slow}.
The external gravity-wave subsystem is given by removing all forcing terms from the momentum equation except for the pressure gradient term, e.g. \eqref{eqn:grav_wave_model}.
Essentially, the splitting can be thought of as advancing the fast tendency terms with a given LTS scheme, and the slow tendency terms with a simple Forward Euler step.
In practice, the splitting works by evaluating the slow tendencies at time \( t = t^n \) and then skipping their update during the normal local time-stepping procedure.
In particular, with a time discretization of the momentum equation in \eqref{eqn:odesys}, the right hand side term \( \Phi \) at a generic time \( t^* \), i.e. \( \Phi\left( \u(t^*),\, h(t^*) \right) \), is approximated as
\begin{linenomath}
\begin{equation}
\begin{aligned}
    \Phi\left( \u(t^*),\, h(t^*) \right) \approx \Phi^{\text{fast}}\left( \u(t^*),\, h(t^*) \right) + \Phi^{\text{slow}}\left( \u(t^n),\, h(t^n) \right) \,.
\end{aligned} \label{eqn:splitting}
\end{equation}
\end{linenomath}
The approximation above can be seen as an additive Lie-Trotter splitting.
Note that because the external gravity-wave system is treated as fast, there is no splitting within the thickness equation, only in the momentum equation.
That is,
\begin{linenomath}
\begin{equation}
\begin{aligned}
    \Phi^\text{fast}(\u,\, h) &= -g\nabla (h + z_{b}) \\
    \Psi^\text{fast}(\u,\, h) &= - \nabla \cdot (h\u)\,.
    \label{eqn:fast_terms}
\end{aligned}
\end{equation}
\end{linenomath}

As a result of this splitting, the more computationally expensive slow terms are only evaluated once per coarse time-step, saving a significant number of floating point operations and MPI communication calls.
This splitting does introduce a first-order error into the temporal local truncation error of our models, but we will show in Section \ref{subsec:hurricane_model} that the quality of the solution is not affected (Figure \ref{fig:split_unsplit}).
This is due to the fact that in complex real-world applications driven by observed data, like the one shown in Section \ref{subsec:hurricane_model}, errors coming from the spatial discretization, model parameterization, and observed data are dominant.
Further, the slow dynamics evolve on the order of hours, so taking time-steps on the order of hundreds of seconds, as in Section \ref{subsec:hurricane_model}, does not introduce meaningful errors.

As mentioned above in Section \ref{subsec:parallalization}, the splitting also helps with the practical concern of needing a way to restrict tendency calculations to certain regions of a mesh for a LTS scheme.
MPAS-Ocean tendency routines are very complex and further rely on other routines that compute diagnostic quantities.
In order to efficiently compute tendencies only on given LTS regions without using multiple memory blocks, one would need to reimplement the MPAS-Ocean tendency routines from scratch in a way that allows the domain of the calculation to be restricted.
This would have been impractical to do for the full tendency routines, but because of the operator splitting, we only need to reimplement the parts of the routines for the fast terms.
The versions of FB-LTS and LTS3 implemented in E3SM are implemented in the context of this operator splitting. 
In the numerical experiments in Section \ref{sec:numerical_experiments}, we use the labels FB-LTS and LTS3 to refer to these schemes with this operator splitting.
In the convergence test from Section \ref{subsec:convergence}, the same implementation of FB-LTS was used, but on a problem with all the slow terms disabled.
This effectively removes the operator splitting, so FB-LTS showes the expected second-order convergence.

A final point of interest is how this splitting effects the CFL performance of the model.
That is, we are interested in how this splitting effects the size of the admittable time-step in both the fine and coarse regions.
In Section \ref{subsec:computational_performance}, we show two cases, one in which the CFL performance is not effected and one in which it is, and provide evidence that this depends on the ratio of the resolution of cells in the coarse region to the resolution of cells in the fine region.


\section{Numerical Experiments}
\label{sec:numerical_experiments}

Here, we present a series of numerical experiments that showcase the performance of FB-LTS as compared to LTS3 and RK4.
As described in Section \ref{subsec:splitting}, both LTS3 and FB-LTS are used in the context of the first-order operator splitting of the fast/slow subsystems within the SWEs; for readability, we continue to refer to the methods as LTS3 and FB-LTS.  
For these experiments, we model the storm surge caused by hurricane Sandy in and around Delaware Bay off the eastern coast of the United States, running on meshes that have been regionally refined to very high resolutions (2 km and 125 m) around Delaware Bay (Table \ref{tbl:meshes} and Figure \ref{fig:meshes_cellWidth}).


\subsection{Hurricane Model}
\label{subsec:hurricane_model}

The momentum and thickness equations for the hurricane Sandy model are given by
\begin{linenomath}
\begin{equation}
\begin{aligned}
    \pd{\u}{t}
        + \left(\nabla \times \mathbf{u} + f \mathbf{k}\right) \times \mathbf{u}  
        =
        & - \nabla K  
          - \frac{1-\beta}{\rho_0}\nabla p^s  
          - g\nabla\left(\xi - \xi_{EQ} - \beta\xi\right) \\[7pt] 
        & - \chi \frac{\mathcal{C} \mathbf{u}}{H}  
          - \mathcal{C}_{\text{D}}\frac{\abs{\mathbf{u}}\mathbf{u}}{h}  
          + \mathcal{C}_{\text{W}}\frac{\abs{\mathbf{u}_{\text{W}} - \mathbf{u} }\left( \mathbf{u}_{\text{W}} - \mathbf{u} \right)}{h}  
          \\[7pt]
    \pd{h}{t} + \nabla \cdot (h \mathbf{u}) &= 0 \,,
    \label{eqn:hurricane_model}
\end{aligned}
\end{equation}
\end{linenomath}
where
$ \mathbf{u} $ is the horizontal fluid velocity,
$ t $ is the time coordinate, 
$ f $ is the Coriolis parameter, 
$ \mathbf{k} $ is the local vertical unit vector, 
$ K = \frac{\abs{\mathbf{u}}^2}{2} $ is the kinetic energy per unit mass,
$ \beta $ is the self-attraction and loading coefficient \citep{accad1978}, 
$ \rho_0 $ is the (constant) fluid density,
$ p^s $ is the surface pressure, 
$ g $ is the gravitational constant, 
$ \xi $ is the sea-surface height perturbation,  
$ \xi_{\text{EQ}} $ is the sea-surface height perturbation due to equilibrium tidal forcing \citep{arbic2018},
$ \frac{\mathcal{C}}{H} $ is a spatially varying internal tide dissipation coefficient \citep{jayne2001}, 
$ \chi $ is a scalar tuning factor optimized for barotropic tides response \citep{barton2022}.
$ H $ is the resting depth of the ocean,
$ h $ is the total ocean thickness such that $ h = H + \xi $,
$ \mathcal{C}_{\text{D}} $ is the bottom drag coefficient, 
$ \mathcal{C}_{\text{W}} $ is the wind stress coefficient \citep{garratt1977},
and $ \mathbf{u}_{\text{W}} $ is the horizontal wind velocity.
Here, the thickness equation is the conservation of volume for an incompressible fluid, where the volume is normalized by the cell area, which is constant in time.

\begin{figure}
    \centering
    \includegraphics[width=0.7\textwidth]{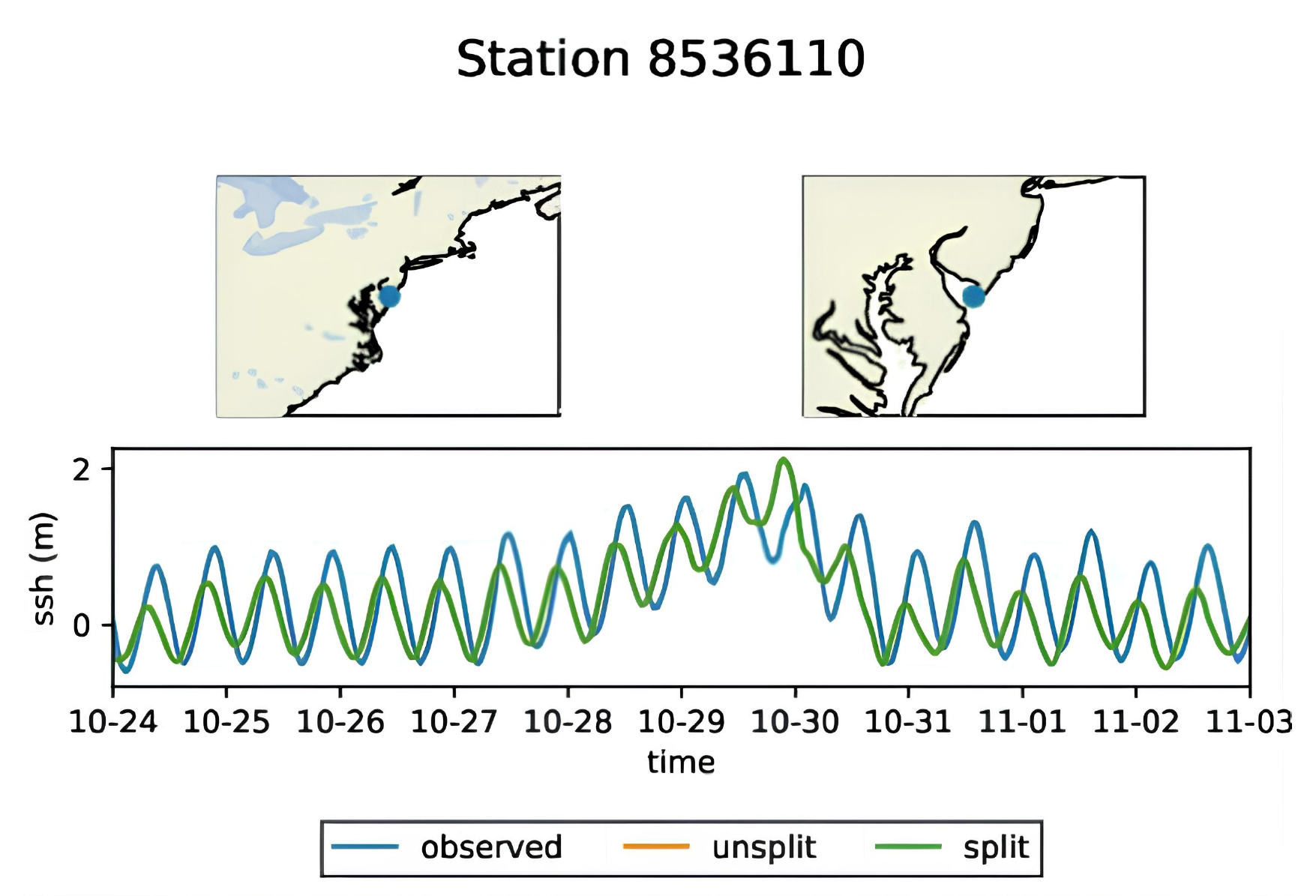}
    \caption{
        The sea-surface height solution produced by the hurricane Sandy model using LTS3 both with and without the operator splitting described in Section \ref{subsec:splitting}.
        The unsplit solution is covered by the split solution, showing the split model produces a solution of the same quality.
    }
    \label{fig:split_unsplit}
\end{figure}

The LTS algorithms implemented in MPAS-Ocean use the operator splitting described in Section \ref{subsec:splitting}; for this hurricane model the terms treated as fast in \eqref{eqn:hurricane_model} are
\begin{linenomath}
\begin{equation}
\begin{aligned}
    \Phi^\text{fast}(\u,\, h) &= -g\nabla\,(\xi - \beta\xi) \\
    \Psi^\text{fast}(\u,\, h) &= - \nabla \cdot (h\u)\,.
    \label{eqn:hurricane_fast_terms}
\end{aligned}
\end{equation}
\end{linenomath}
To show that this splitting does not degrade the quality of the model solution, we run the hurricane model using LTS3 from the now outdated version of MPAS-Ocean used in \citet{lilly2023_sandy}, and the current version (Figure \ref{fig:split_unsplit}).
Again, when we refer to FB-LTS or LTS3 in this section, we specifically mean a time-stepping scheme using the operator splitting from Section \ref{subsec:splitting} with the given LTS scheme.

\begin{figure}
    \centering
    \begin{tabular}{lcc}
        & \multicolumn{2}{c}{\includegraphics[width=0.9\textwidth]{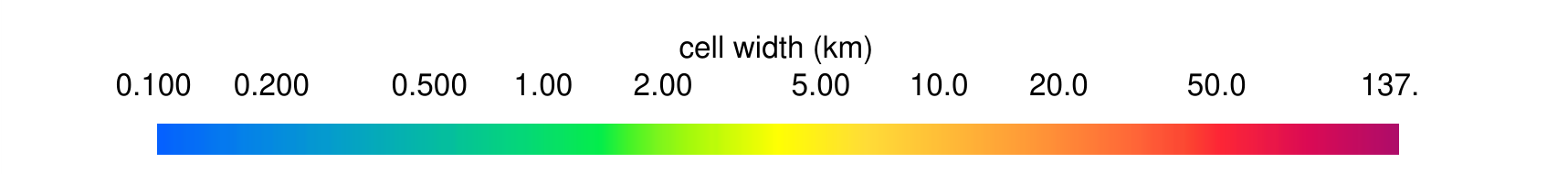}} \medskip \\
        & \large Atlantic Ocean & \large Eastern US Coast \smallskip \\
        \rotatebox{90}{\hspace{0.3cm} \large DelBay2km}
            & \includegraphics[width=0.44\textwidth]{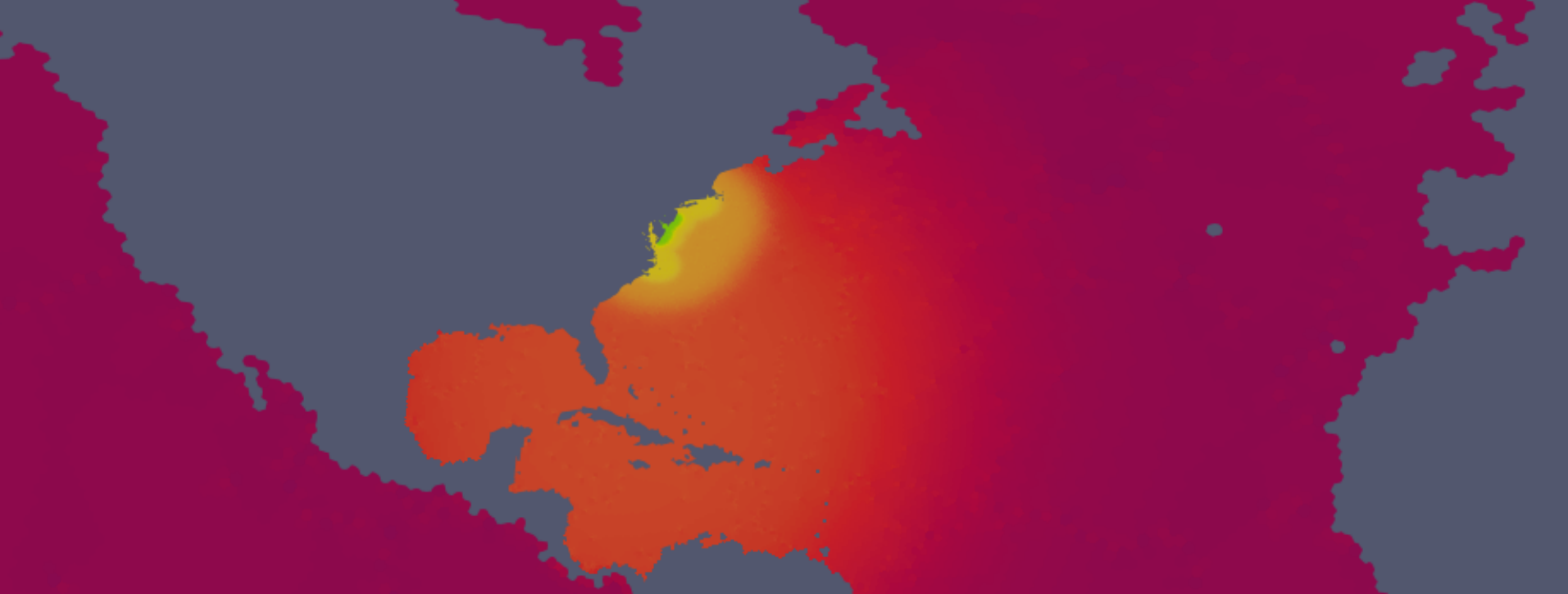}
            & \includegraphics[width=0.44\textwidth]{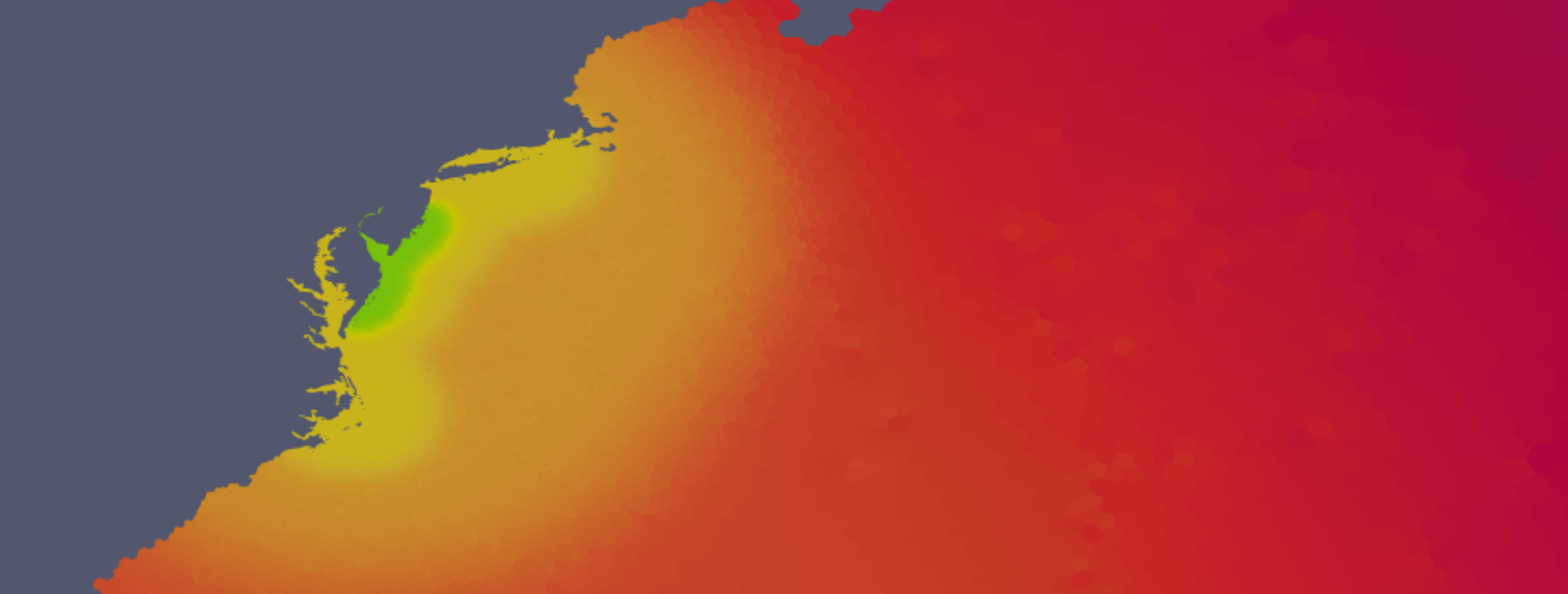} \medskip \\
        \rotatebox{90}{\hspace{0.2cm} \large DelBay125m}
            & \includegraphics[width=0.44\textwidth]{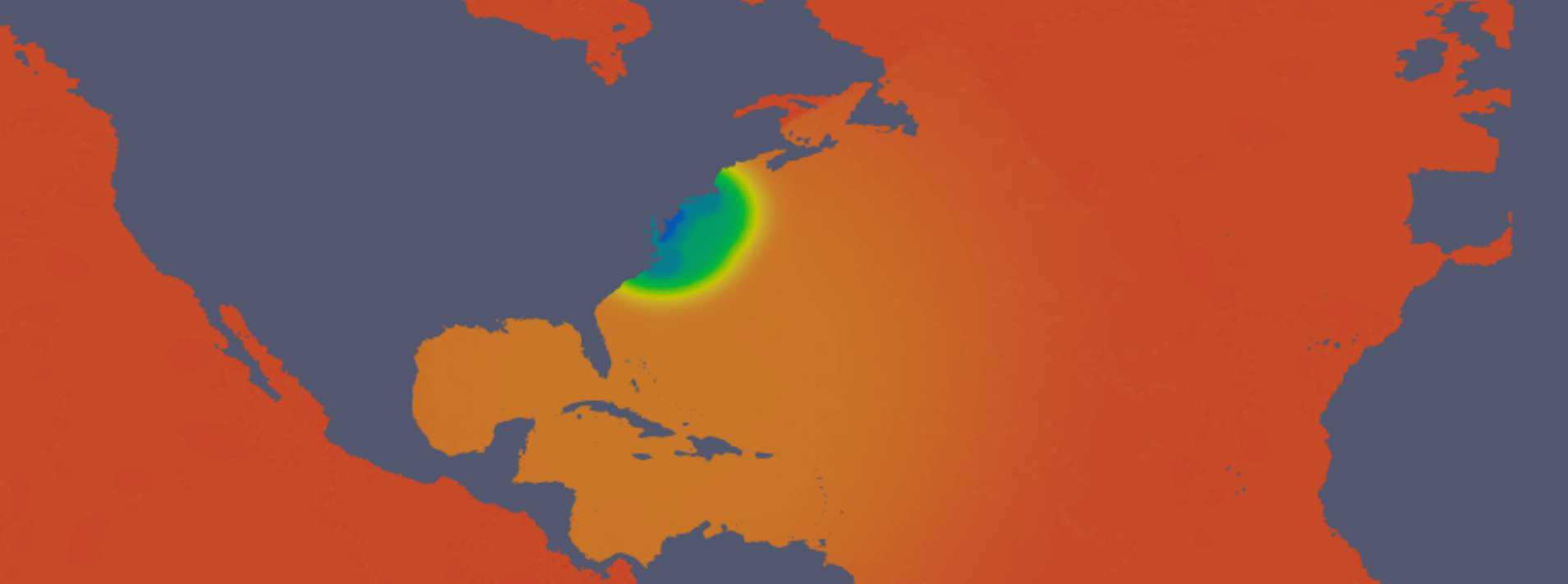}
            & \includegraphics[width=0.44\textwidth]{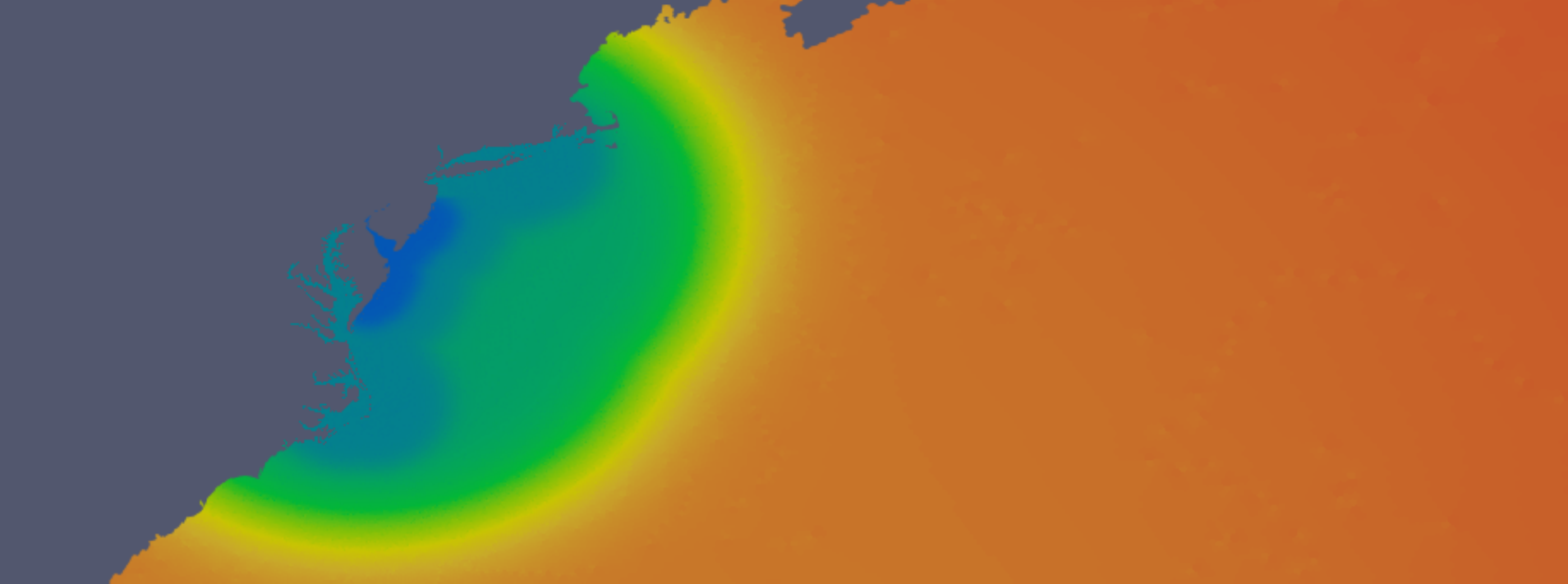} \\
    \end{tabular}
    
    \caption{
        Cell width in kilometers for DelBay2km and DelBay125m as described in Table \ref{tbl:meshes}.
        Each subplot shares the same color scale in log-space.
        Note that these are global meshes and that parts of the globe not pictured here use the background resolutions shown in Table \ref{tbl:meshes}.
    }
    \label{fig:meshes_cellWidth}
\end{figure}

We consider two meshes, DelBay2km and DelBay125m, shown in Figure \ref{fig:meshes_cellWidth} and described in Table \ref{tbl:meshes}.
These meshes are similar to those by the same names from \citet{lilly2023_sandy}, but have slight differences due to updates to mesh generation libraries.
Both meshes highly resolve Delaware Bay, then smoothly vary in resolution out to a low global background resolution.
These meshes were chosen to highly resolve the hurricane itself while controlling the overall number of cells in the global mesh, making the potential benefits of local time-stepping schemes obvious.
The simulations are run for 24 simulated-days; starting on 10/10/2012 and ending on 11/03/2012.

Both LTS3 and FB-LTS were configured to use the same configuration of the LTS regions, taking the fine region (the region where the fine time-step is used) to be around the Delaware coast and the coarse region (the region where the coarse time-step is used) to be the rest of the globe as pictured in Figure \ref{fig:meshes_cellWidth}.
One can note that FB-LTS scheme uses more cells labeled as interface adjacent fine cells than LTS3; these are the cells denoted as \( \F_P^5 \) cells in Section \ref{subsec:fblts}.
These cells are advanced with the same fine time-step as the rest of the fine region, they are only differentiated within the FB-LTS and LTS3 algorithms because these cells require a small number of additional computations in order to obtain needed data on the interface one and two regions during the coarse advancement step of both methods.
In order to obtain the data needed to perform the interface one prediction step, we need to obtain uncorrected data on interface one at time \( t^{n+1} \).
FB-RK(3,2) takes three stages to obtain data at time \( t^{n+1} \), while SSPRK3 produces a prediction for this data with its first stage; as a result the domain of dependence for interface one cells to obtain this needed data is larger for FB-LTS.
One can also note that, as described in Section \ref{subsec:fblts}, the FB-LTS algorithm requires these extra computations on a decreasing subset of these interface adjacent fine cells as the method progresses through the coarse advancement step, starting with \( \F_P^5 \), then \( \F_P^4 \), then \( \F_P^3 \), et cetera. 
In order to decrease the complexity of the FB-LTS implementation in MPAS-Ocean we opt not to shrink this region and instead perform these extra calculations on all the interface adjacent cells pictured in Figure \ref{fig:lts_regions}.
This has no effect on the numerical scheme itself as the unneeded data is simply thrown away, and the penalty to the computational performance of the implementation is negligible.

\begin{figure}
    \centering

    \includegraphics[width=0.95\textwidth]{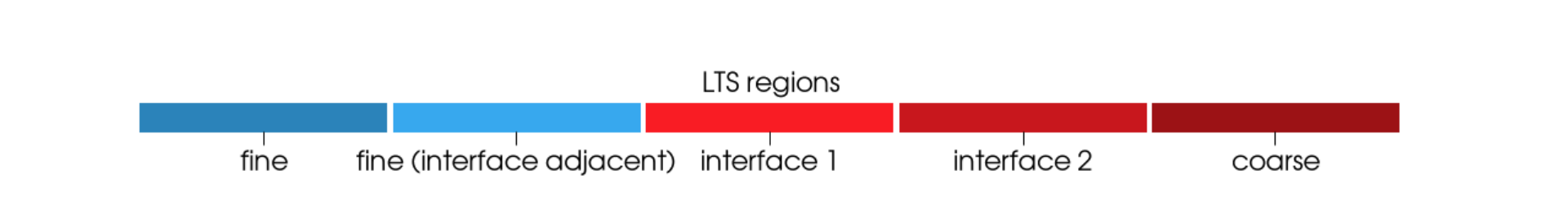}
    \vspace{-5pt}
    
    \begin{subfigure}{0.49\textwidth}
        \includegraphics[width=0.9\textwidth]{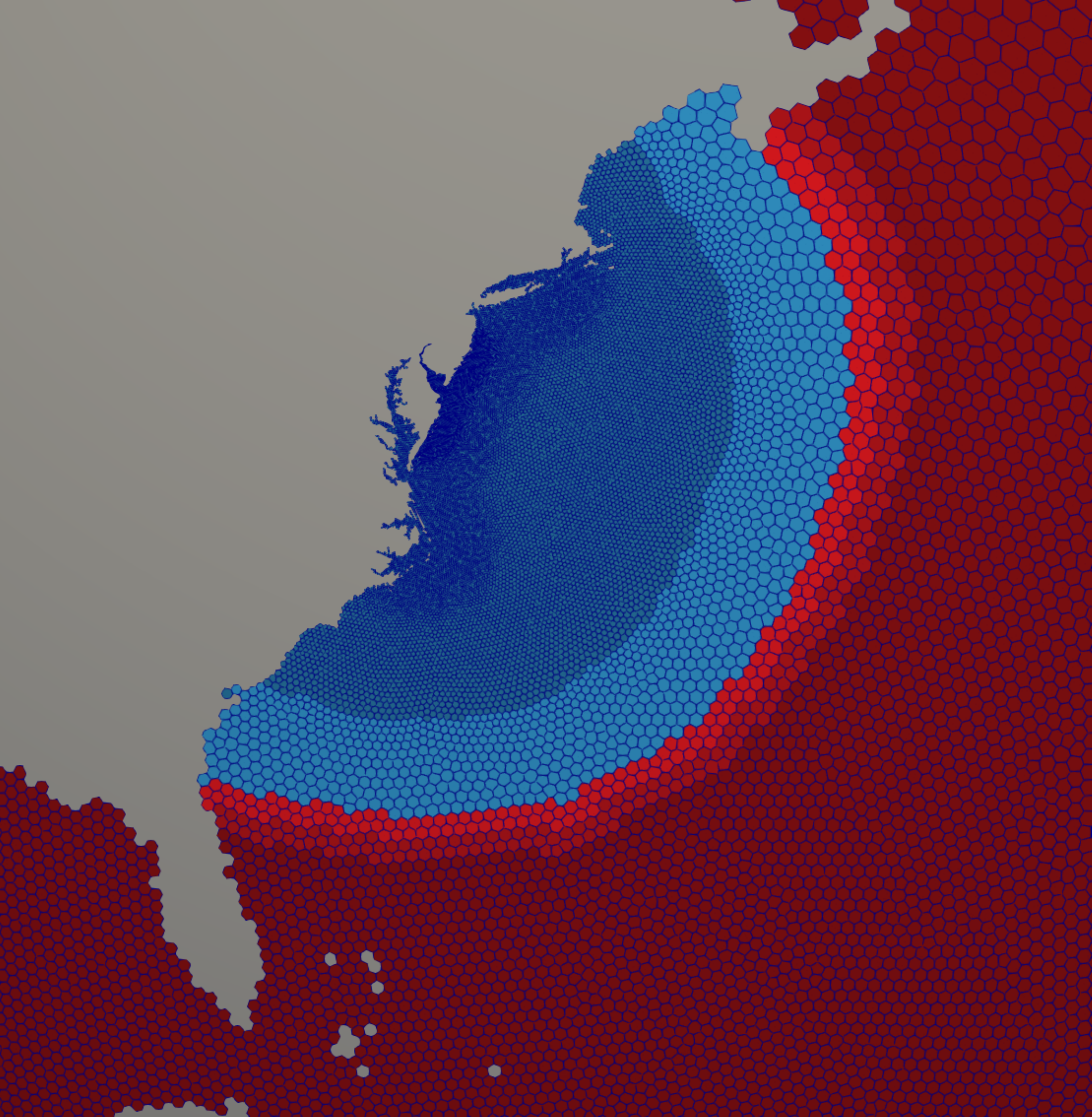}
        \caption{LTS regions for FB-LTS}
        \label{fig:fblts_regions}
    \end{subfigure}
    \begin{subfigure}{0.49\textwidth}
        \includegraphics[width=0.9\textwidth]{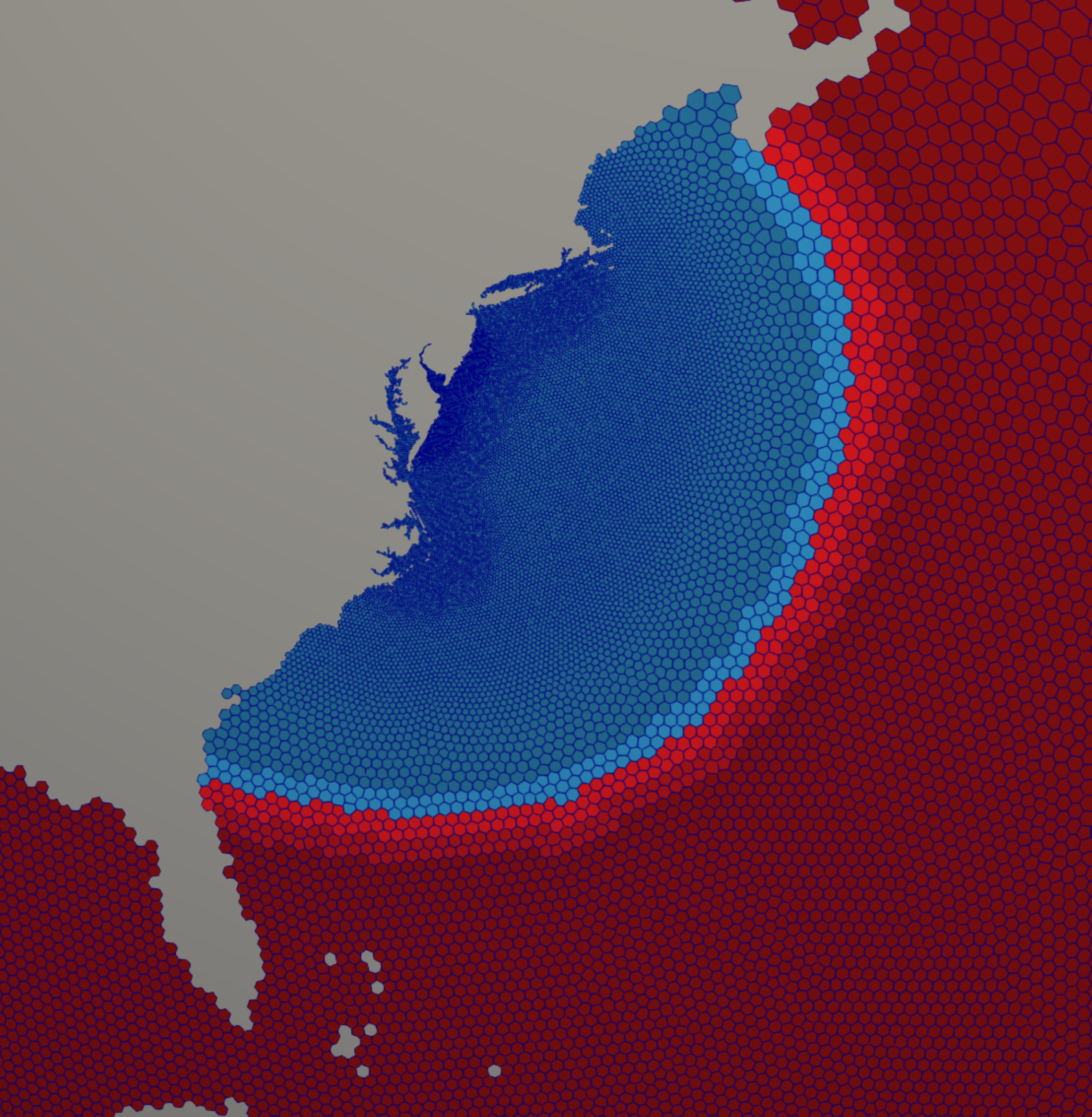}
        \caption{LTS regions for LTS3}
        \label{fig:lts3_regions}
    \end{subfigure}
    \caption{
        Configurations for the LTS regions for FB-LTS and LTS3 as shown on DelBay2km.
        The LTS regions are configured similarly on DelBay125m.
    }
\label{fig:fblts_lts3_regions}
\end{figure}

The time-steps used on each mesh by each LTS scheme for performance tests are given in Table \ref{tbl:meshes}.
These time-steps were obtained experimentally by running the model at increasing time-steps until it becomes unstable and selecting the largest time-steps for which it is stable.
In the case of DelBay2km, the model is stable at the reported time-steps for the entire 24 simulated-day duration of the hurricane simulation.
Because of the greatly increased computational cost of running the model on DelBay125m, the reported time-steps are obtained by running for one simulated-day starting from an initial condition given by the model state after three simulated days (10/13), which is when the tidal dynamics stabilize after spinning up from rest.
When finding the fine and coarse time-steps in this way, there are two possible approaches; the fine and coarse time-steps must differ by an integer factor of \( M \) such that \( \Delta t_{\text{fine}} = \frac{\Delta t_{\text{coarse}}}{M} \), so they can be obtained in different orders.
One could find the largest time-step admittable on the coarse region of the mesh then find the smallest value of \( M \) that gives an admittable fine time-step, or first find the largest admittable fine time-step then the largest value of \( M \) that gives an admittable coarse time-step.
This has the result of the time-step in one region or the other not technically being maximal.
For both DelBay2km and DelBay125m, we opt to first maximize the time-step in the fine region, then find the largest admittable value of \( M \).

\begin{table}
    \centering
        \begin{tabular}{lrr}
            ~ & \textbf{DelBay2km}  & \textbf{DelBay125m} \\ \toprule
            \textbf{Grid Cell Width} (km)  & ~ & ~ \\ \midrule
            Global background & 120 & 30 \\ 
            Western Atlantic & 30 & 15  \\ 
            Eastern US coast & 10 & 0.625 \\ 
            Delaware coast & 5 & 0.3125  \\ 
            Delaware Bay & 2 & 0.125  \\ \toprule
            \textbf{Mesh Parameters} & ~ & ~  \\ \midrule
            Number of cells & 58,141 & 4,617,372  \\
            Count ratio  & 1.92 & 0.12  \\
            Resolution ratio & 15 & 120  \\ \toprule
            \textbf{Time-Steps} & ~ & ~ \\ \midrule
            RK4 & ~ & ~ \\
            \( \Delta t_{\text{global}} \) (s) & 34 & 2  \\ \midrule
            LTS3 & ~ & ~ \\
            \( \Delta t_{\text{fine}} \) (s) & 20 & 1  \\ 
            \( \Delta t_{\text{coarse}} \) (s) & 60 & 4  \\
            \( M \) & 3 & 4  \\ \midrule
            FB-LTS & ~ & ~ \\
            \( \Delta t_{\text{fine}} \) (s) & 46 & 2  \\ 
            \( \Delta t_{\text{coarse}} \) (s) & 138 & 4  \\
            \( M \) & 3 & 2  \\ \midrule
            Time-Step Ratios & ~ & ~ \\
            \( \nicefrac{ \Delta t_\text{fine}^\text{FB-LTS} }{ \Delta t_\text{global}^\text{RK4} } \) & 1.35 & 1.00 \\
            \( \nicefrac{ \Delta t_\text{fine}^\text{FB-LTS} }{ \Delta t_\text{fine}^\text{LTS3} } \) & 2.30 & 2.00 \\
            \( \nicefrac{ \Delta t_\text{coarse}^\text{FB-LTS} }{ \Delta t_\text{coarse}^\text{LTS3} } \) & 2.30 & 1.00 \\ \bottomrule
        \end{tabular}
    \caption{
        Relevant parameters for each mesh and LTS scheme used in performance experiments.
        The \textit{count ratio} is the ratio of the number of coarse cells to the number of fine cells in the mesh, i.e. a count ratio greater than one means that the mesh contains more cells using the coarse time-step than the fine time-step.
        The \textit{resolution ratio} is the ratio of the cell width of the coarse cells to the cell width of the fine cells. 
        In the case where there are cells of multiple resolutions in either region, as is the case in our meshes, we consider the smallest value of cell width in a given region as it is the smallest cell that restricts the size of the time-step admittable in that region.
        For both meshes, the fine region includes Delaware Bay, the Delaware coast, and the Eastern US coast.
        The coarse region includes the Western Atlantic and the rest of the globe (global background). 
    }
    \label{tbl:meshes}
\end{table}


\subsection{Computational Performance}
\label{subsec:computational_performance}

In order to demonstrate the computational performance of FB-LTS, we run the hurricane Sandy model using RK4, LTS3, and FB-LTS on both DelBay2km and DelBay125m.
For each mesh, we run the model for a number of simulated-seconds equal to the least common multiple of the time-steps used by each scheme.
We measure the speedup by taking the ratio of the CPU-time for the slower scheme to that of FB-LTS, that is,
\begin{linenomath}
\begin{equation}
\begin{aligned}
     \text{speedup vs. RK4} &= \frac{\text{RK4 CPU-time}}{\text{FB-LTS CPU-time}} \\
     &\null \\
     \text{speedup vs. LTS3} &= \frac{\text{LTS3 CPU-time}}{\text{FB-LTS CPU-time}} \,.
    \label{eqn:speedup}
\end{aligned}
\end{equation}
\end{linenomath}
On DelBay2km, FB-LTS outperforms RK4 by a factor of 10.08 and outperforms LTS3 by a factor of 2.27, and on DelBay125m FB-LTS outperforms RK4 by a factor of 5.13 and outperforms LTS3 by a factor of 1.32 (Table \ref{tbl:performance}).
Previous to the implementation of LTS3 and FB-LTS in the MPAS-Ocean code-base, RK4 scheme was the method of choice for single-layer configurations like that used in our hurricane model, so the speedup versus RK4 achieved by both LTS schemes is significant and points to the value of local time-stepping schemes and the operator splitting described in Section \ref{subsec:splitting} for efficiency of ocean simulations.
FB-LTS further improves on the speedup achieved by LTS3 by taking advantage of the CFL optimized scheme on which it is based, and it will be shown in Section \ref{subsec:solution_quality} that this additional speedup does not incur a penalty to the quality of the model solution.

\begin{table}
    \centering
    \begin{tabular}{lrr}
        ~ & \textbf{DelBay2km} & \textbf{DelBay125m} \\ \toprule
        run-time (hh:mm:ss) & 06:31:00 & 00:10:00 \\ 
        number of MPI ranks & 32 & 128 \\ \midrule
        RK4 (s) & \textit{26.07} & \textit{186.02} \\ \bottomrule
        LTS3 (s) & \textit{5.86} & \textit{47.85} \\
        speedup vs. RK4 & 4.45 & 3.89 \\ \bottomrule
        FB-LTS (s) & \textit{2.59} & \textit{36.26} \\ 
        speedup vs. RK4 & 10.08 & 5.13 \\
        speedup vs. LTS3 & 2.27 & 1.32
    \end{tabular}
    \caption{
        CPU-time performance of RK4, LTS3, and FB-LTS on DelBay2km and DelBay125m.
        Each reported time above is the average of individual runs, to account for conditions on the machine.
        The run-times were chosen as common multiples of the time-steps used be each scheme as reported in Table \ref{tbl:meshes}.
    }
    \label{tbl:performance}
\end{table}

We can also compare the results from Table \ref{tbl:performance} to results from Table 2 from \citet{lilly2023_sandy}.
On DelBay2km, the previous unsplit implementation of LTS3 only outperformed RK4 by a factor of 1.48, while on DelBay125m unsplit LTS3 was slower than RK4 by a factor of 1.16.
One should keep in mind that there multiple factors involved in the difference in performance between these two cases; the new schemes benefit from less MPI communication overhead due to the improved domain decomposition paradigm described in Section \ref{subsec:parallalization} and the operator splitting described in Section \ref{subsec:splitting} in particular.
Regardless, this shows the strong progression of efforts to increase the computational efficiency of these SWE solvers.

As noted in Section \ref{subsec:splitting}, we are also interested in how the operator splitting affects both the computational performance and the CFL performance of our LTS schemes.
On DelBay2km, the results suggest that the CFL performance of both LTS3 and FB-LTS is not affected by the splitting.
That is, the restriction placed on the maximum admittable time-steps, as reported in Table \ref{tbl:meshes}, are enforced by the terms treated as fast.
In this case FB-LTS is taking time-steps 2.3 times larger than LTS3 in both the fine and coarse regions, as we expect.
Contrast this with the time-steps used on DelBay125m; FB-LTS outperforms LTS3 in CFL efficiency by a factor of 2 in the fine region but both take a time-step of 4 s in the coarse region.
This behavior is not consistent with what we expect to achieve assuming the time-steps are being bound by the fast terms as we would hope.
Further, we would expect that the ratio \( M \) between the fine and coarse time-step was much larger.
For DelBay2km the resolution ratio is 15 and we get \( M = 3 \).
On DelBay125m with a resolution ratio of 120 and the same underlying model equations, we would expect to have the \( M \) is at least 3, if not significantly larger (in \citet{lilly2023_sandy}, the authors found that \( M = 24 \) for LTS3 in this case).

Motivated by the desire to understand this unexpected result, we ran additional tests that revealed that the problem was due to the global CFL restriction imposed by the slow terms.
Running the model on DelBay125m, with both LTS3 and FB-LTS, we observed that any increase to the coarse time-step would cause the model to become unstable, with the instability manifesting in the fine region far from the interface regions, deep in Delaware Bay.
The fact that increasing the coarse time-step causes instability in the fine region means that it is very unlikely that the instability is being caused by the fast terms, which are being advanced with LTS.
The way the slow terms are advanced by the split scheme can be thought of as advancing the slow part of the system by a single Forward Euler step using the coarse time-step everywhere on the mesh.
This explains why increasing the coarse time-step causes instability in the fine region, where the cells are very small and the CFL condition therefore very restrictive.
Effectively, the operator splitting places a bound on the size of the coarse time-step that depends on the resolution of the fine cells, possibly limiting the effectiveness of the LTS schemes employed in this splitting.
As the results from Table \ref{tbl:performance} show, this does not mean that the split algorithm using LTS cannot be computationally efficient, it only means that it can limit the CFL performance of LTS schemes in cases of extreme differences in resolution, like that in DelBay125m.

A possible way to combat this limit on CFL performance imposed by the splitting would be to alter the algorithm so that, during the \textit{Fine Advancement} step, the slow tendencies were evaluated once per fine-step. 
Instead of using the values of the slow tendencies as evaluated at time \( t^n \) for the entire routine, we could calculate the slow tendencies at times \( t^{n,k} \) for \( k = 0, \cdots M-1 \) and use these values while advancing on fine cells.
Doing this, the algorithm would effectively be advancing the slow subsystem still with forward Euler, but with the fine time-step rather than the coarse time-step.
Exploration of this particular change to the operator splitting algorithm is outside the scope of this work, as its implementation in MPAS-Ocean would be particularly difficult due to the limitations of the framework described in Section \ref{sec:implementation}.

\begin{figure}
    \centering
    \includegraphics[width=0.7\textwidth]{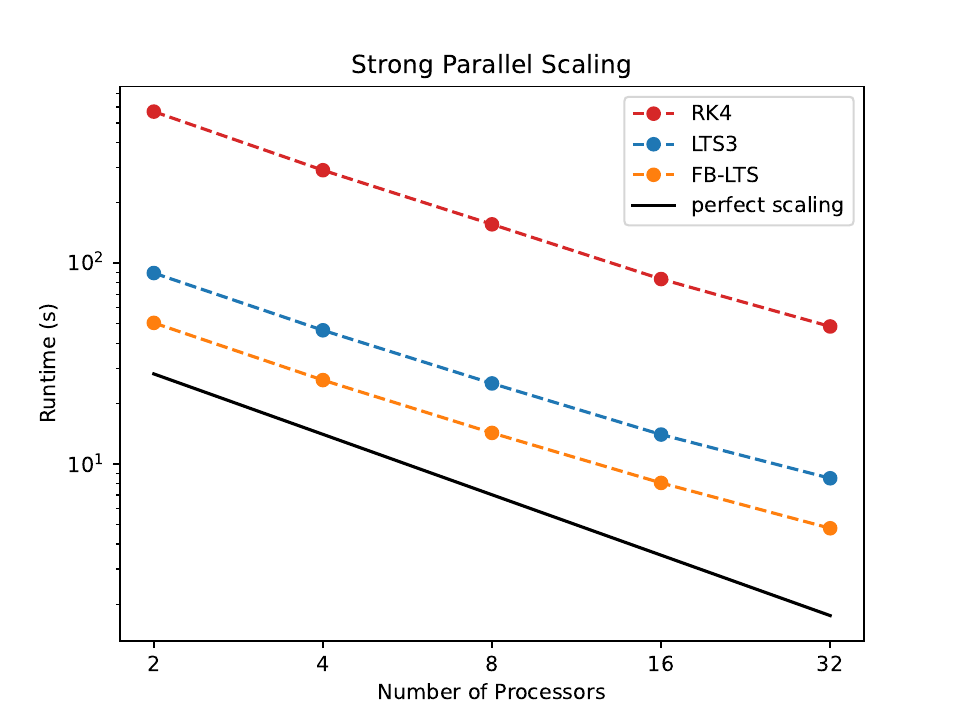}
    \caption{
        Strong parallel scaling in the hurricane Sandy test case on DelBay2km, run for 12 simulated hours.
    }
    \label{fig:parallel_scaling}
\end{figure}

Testing the strong parallel scaling of RK4, LTS3, and FB-LTS, we see the expected approximately linear scaling as the number of MPI ranks is increased (Table \ref{fig:parallel_scaling}).
For each of the three methods, one can notice a slight degradation in the scaling as we reach 32 ranks; this occurs be cause the number of cells per processor (which in this case is approximately 1,800) is small enough that the communication between MPI ranks begins to dominate the simulation.
This is expected behavior, and simply a signal that running the simulation on this mesh with more ranks would have diminishing returns. 

\subsection{Solution Quality}
\label{subsec:solution_quality}

\begin{figure}
    \centering

    \begin{subfigure}{0.49\textwidth}
        \includegraphics[width=\textwidth]{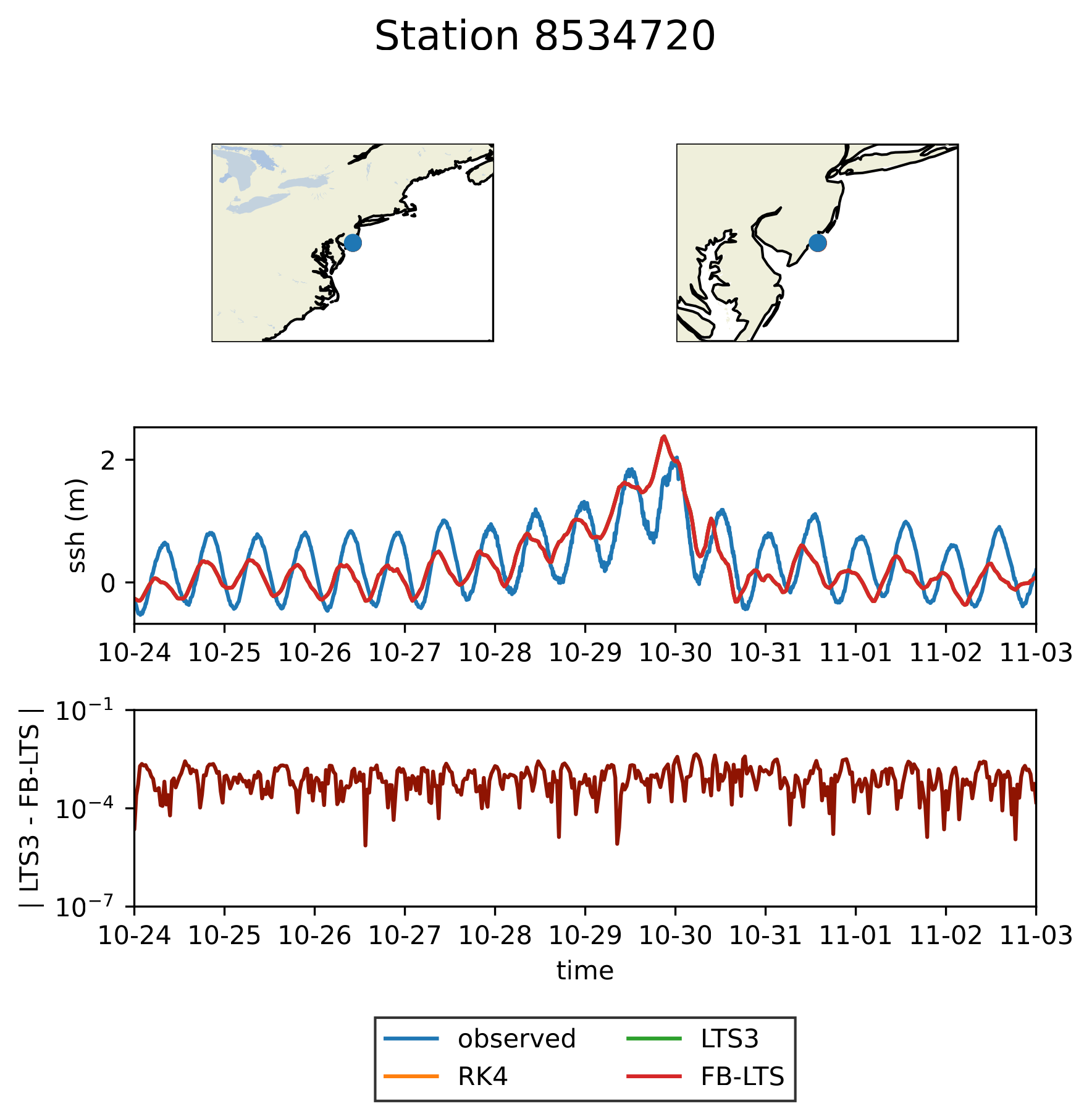}
        \caption{~}
        \label{subfig:8534720}
    \end{subfigure} \hfill
    \begin{subfigure}{0.49\textwidth}
        \includegraphics[width=\textwidth]{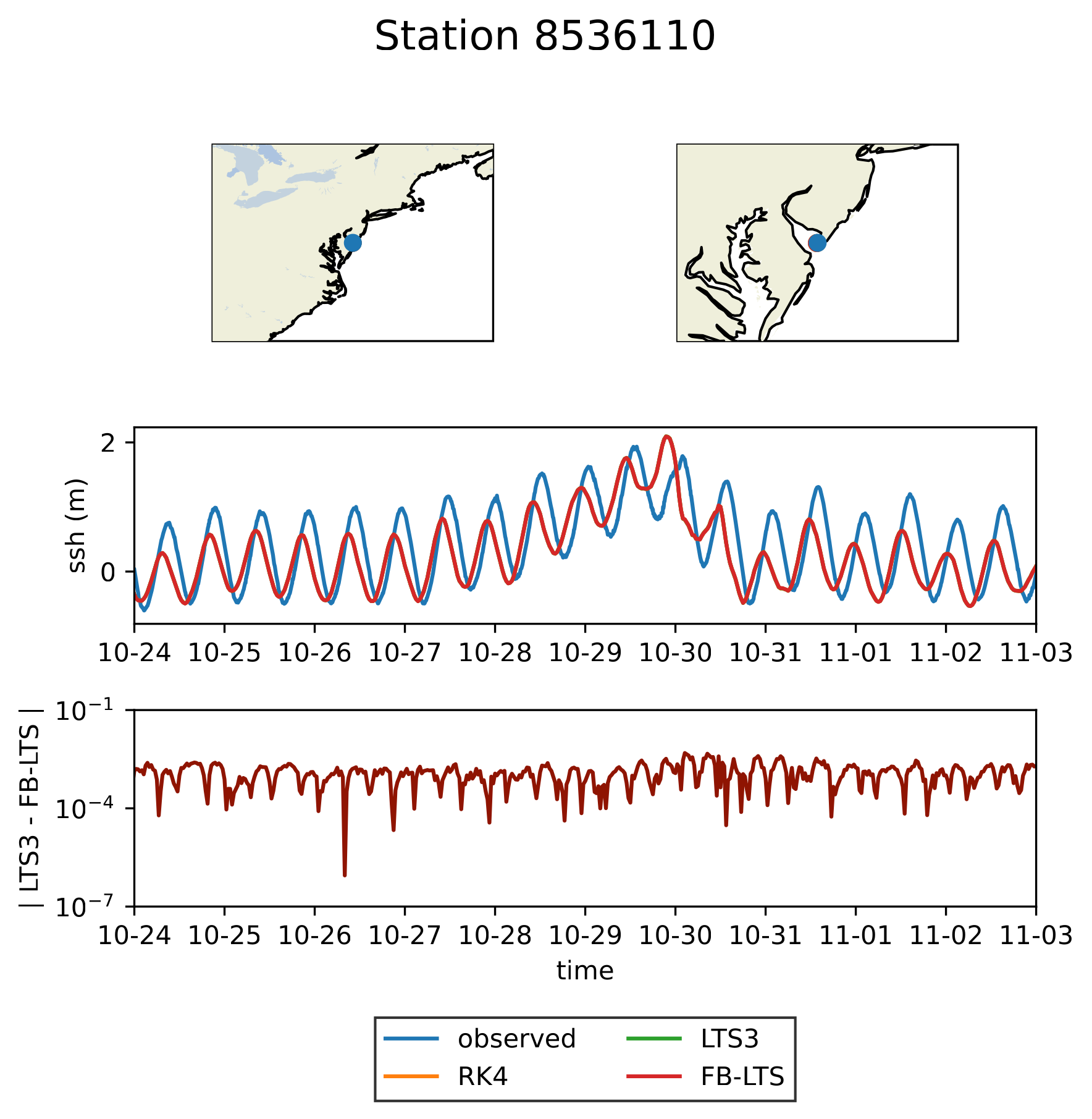}
        \caption{~}
        \label{subfig:8536110}
    \end{subfigure} \\ \bigskip

    \begin{subfigure}{0.49\textwidth}
        \includegraphics[width=\textwidth]{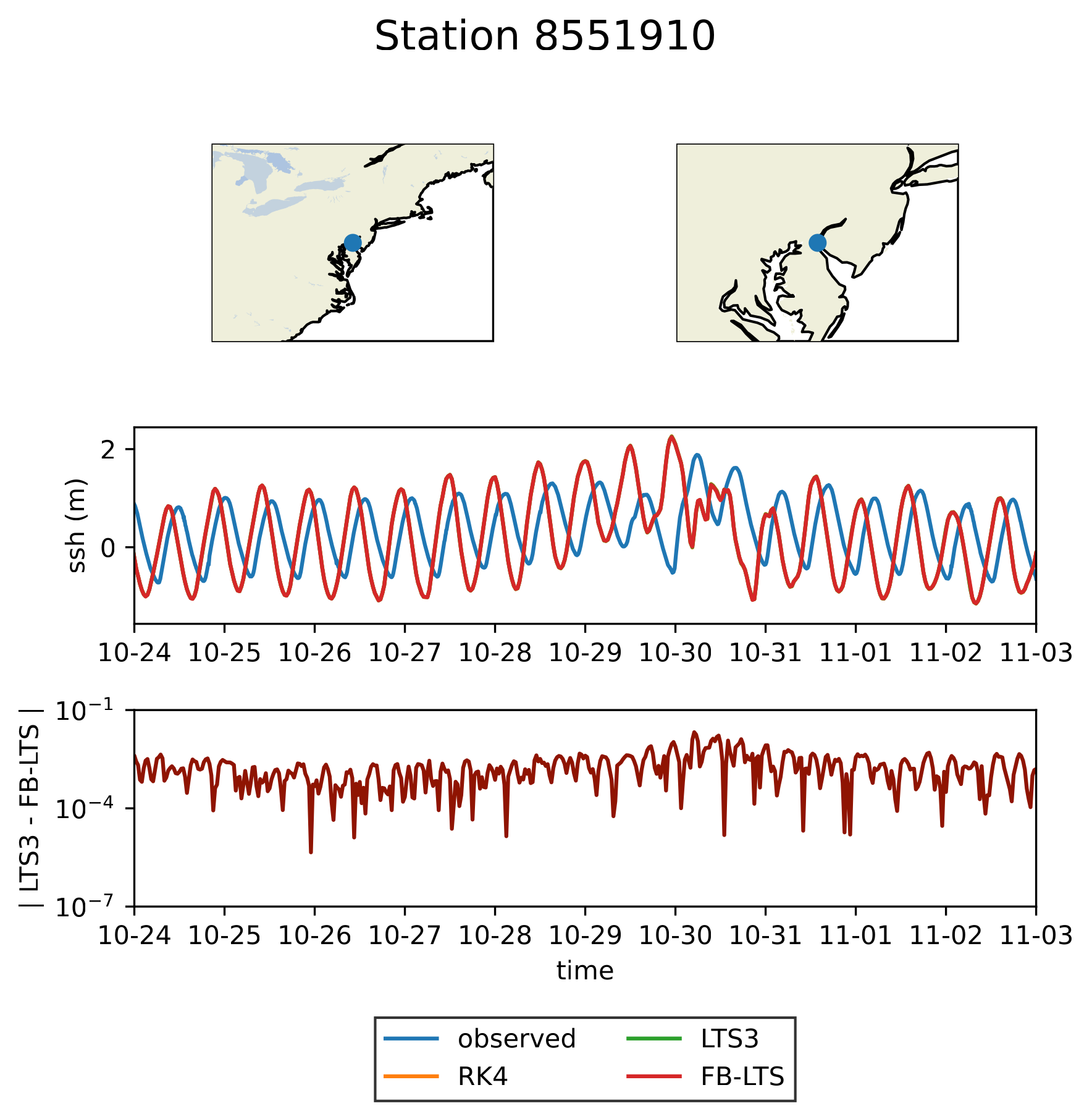}
        \caption{~}
        \label{subfig:8551910}
    \end{subfigure} \hfill
    \begin{subfigure}{0.49\textwidth}
        \includegraphics[width=\textwidth]{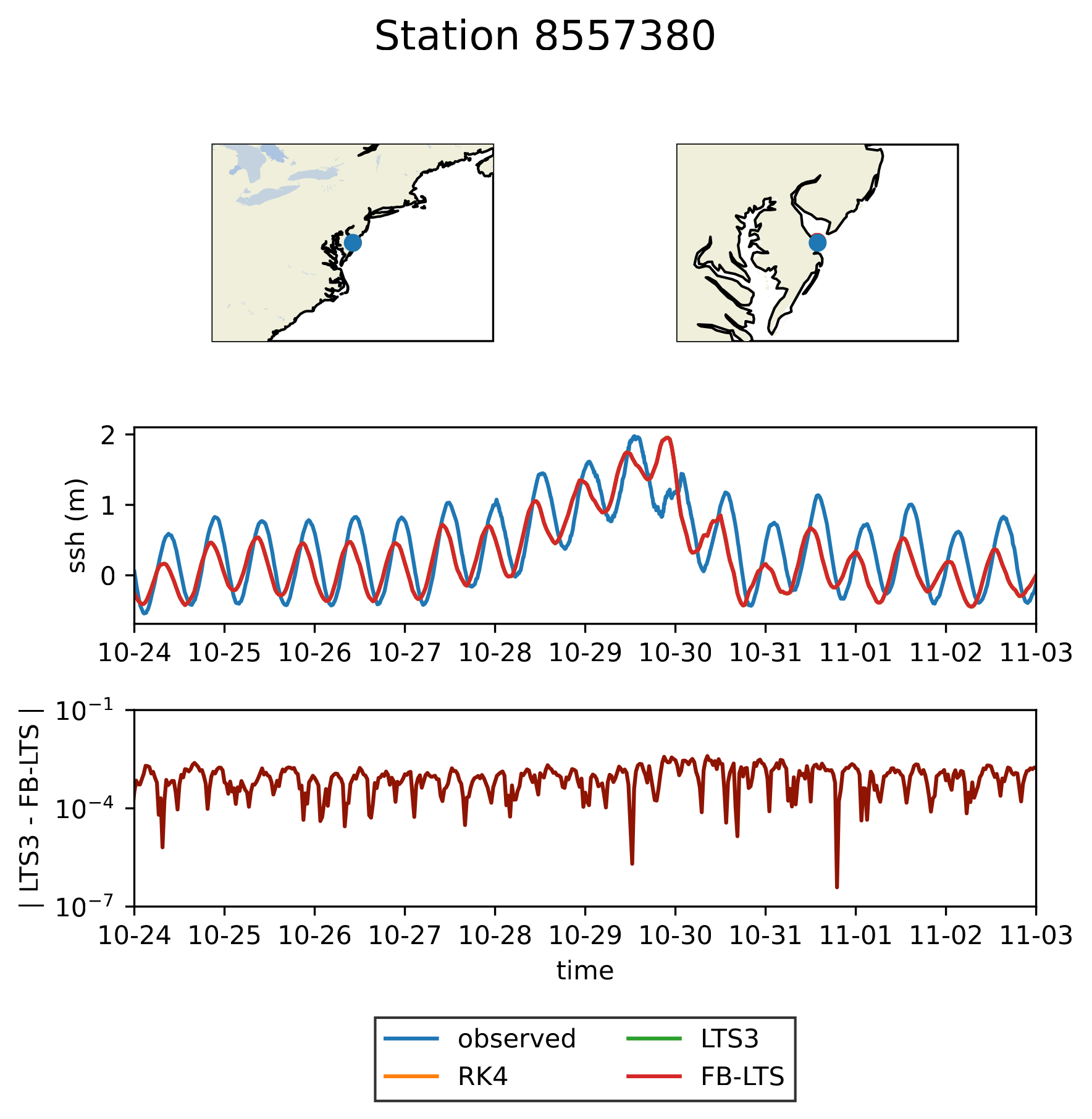}
        \caption{~}
        \label{subfig:8557380}
    \end{subfigure}
    
    \caption{
        SSH solutions from RK4, LTS3, and FB-LTS compared to observed tide gauge data on DelBay2km, using time steps of 30 s, 60 s, and 120 s respectively, with \( M = 3 \) for both LTS schemes.
        Note that The solution curves are directly on top of each other.
        The absolute difference between LTS3 and FB-LTS is shown on a log scale in the lower plots.
    }
    \label{fig:ssh}
\end{figure}

Along with the greatly increased performance discussed in Section \ref{subsec:computational_performance}, FB-LTS is able to produce SSH solutions of the same quality of both RK4 and LTS3.
We run the hurricane model using all three time-stepping schemes on DelBay2km for the full duration of the simulation and record the model SSH at the locations of tidal gauges in and around Delaware Bay.
We also compare model solutions to observed data.
The observed data are from NOAA's Center for Operational Oceanographic Products and Services (CO-OPS) gauges and are available at \url{https://tidesandcurrents.noaa.gov/}.
We observe that the model solutions for each time-stepping scheme do not differ meaningfully.
In particular, the difference between the SSH solutions produced by LTS3 and FB-LTS is, at most, on the order of centimeters (Figure \ref{fig:ssh}).

\section{Conclusion}
\label{sec:conclusion}

We have presented a new local time-stepping scheme (FB-LTS) for the shallow water equations based on the CFL optimized forward-backward Runge-Kutta schemes from \citet{lilly2023_cfl}.
In a real-world test case, FB-LTS produces solutions qualitatively equivalent to those produced by a SSPRK3 based local time-stepping scheme (LTS3) and by the classical four-stage, fourth-order Runge-Kutta method (RK4).
Further, these solutions are produced at a significantly reduced computational cost; the MPAS-Ocean implementation of FB-LTS is up to 10.08 times faster (in terms of CPU-time) than RK4, and up to 2.27 times faster than LTS3.
Numerical experiments show that FB-LTS is a second order scheme everywhere, including on interface cells that allow communication between the fine and coarse regions.
We have also shown that FB-LTS gives exact conservation of mass and absolute vorticity when applied to a TRiSK discretization.

Moving forward, we are interested in adapting FB-LTS for use in a multi-layer ocean model that uses a barotropic-baroclinic splitting.
In particular, FB-LTS would be a strong choice for a barotropic solver, as the barotropic subsystem is similar to the SWEs, so we would expect to see performance benefits very similar to those reported in this work.
With the vital importance of large scale models of the Earth's climate, a solver like FB-LTS would be capable of running complex climate-scale simulations at higher resolutions much faster.



\clearpage

\section*{Acknowledgments}

JRL was supported in part by the National Science Foundation (NSF)
Mathematical Sciences Graduate Internship (MSGI).
GC, MRP, and DE were supported by the Earth System Model Development program area of the U.S. DOE, Office of Science, Office of Biological and Environmental Research as part of the multi-program, collaborative Integrated Coastal Modeling (ICoM) project.


\section*{Data Statement}

The source code for E3SM (and MPAS-Ocean as a component) that implements FB-LTS is available on GitHub.
\begin{itemize}[leftmargin=2.5cm]
    \item[E3SM:] \url{https://github.com/E3SM-Project/E3SM/tree/6b38b86}
\end{itemize}
Similarly, the test case configuration libraries that implement the hurricane Sandy test case (Compass) and the external gravity wave test case (Polaris) for FB-LTS are also available on GitHub.
\begin{itemize}[leftmargin=2.5cm]
    \item[Compass:] \url{https://github.com/MPAS-Dev/compass/tree/74edf14}
    \item[Polaris:] \url{https://github.com/jeremy-lilly/polaris} \\
                    \url{/tree/aquaplanet-external-gravity-wave}
\end{itemize}


\appendix



\section{Derivation of Interface Prediction Coefficients}
\label{sec:derivation_pred_coeffs}

FB-RK(3,2) is a second order scheme, so we will derive second order predictor coefficients for use in the \textit{Interface Prediction} step of FB-LTS.
Assume that we already have data for \( \h^n \) and \( \u^n \), and uncorrected data \( \tilde{\h}^{n+\nicefrac{1}{3}} \), \( \tilde{\u}^{n+\nicefrac{1}{3}} \), \( \tilde{\h}^{n+\nicefrac{1}{2}} \), \( \tilde{\u}^{n+\nicefrac{1}{2}} \), \( \tilde{\h}^{n+1} \), and \( \tilde{\u}^{n+1} \) on interface one.

Start with the spatially discretized system in \eqref{eqn:semi_discrete_sys}, and assume that both \( \Phi_e \) and \( \Psi_i \) are Lipschitz for all \( e \in \C_E^\ifone \) and \( i \in \C_P^\ifone \).
First, we need predictions for \( \h^{n,k} \) for \( k = 0, \cdots M \) and \( \u^{n,k} \) for \( k = 0, \cdots M-1 \).
Start with \( \h^{n,k} \); take a Taylor series expansion centered at time \( t = t^n \), writing terms out to get a second order approximation for \( i \in \C_P^\ifone \),
\begin{linenomath}
\begin{align*}
    h_i\left( t^n + k\frac{\Delta t}{M} \right) &= h_i(t^n) + \left( t^n + k\frac{\Delta t}{M} - t^n \right) \pd{h_i}{t}(t^n) + \Ord{2} \\
    &=  h_i(t^n) + k\frac{\Delta t}{M} \pd{h_i}{t}(t^n) + \Ord{2} \,.
\end{align*}
\end{linenomath}
Next, we approximate \( h_i(t^n) \) by \( h_i^n \), which is a second order approximation (the order of FB-RK(3,2)).
Then \( \pd{h_i}{t}(t^n) \) can be approximated by Forward Euler (FE), \(  \pd{h_i}{t}(t^n) = \frac{\tilde{h}_i^{n+1} - h^n}{\Delta t} + \Ord{1} \).
Insert these into the above,
\begin{linenomath}
\begin{align*}
     h_i\left( t^n + k\frac{\Delta t}{M} \right) &= \left( h_i^n + \Ord{2} \right) + k\frac{\Delta t}{M} \left( \frac{\tilde{h}_i^{n+1} - h^n}{\Delta t} + \Ord{1} \right) + \Ord{2} \\
     &= h_i^n +  k\frac{\Delta t}{M} \frac{\tilde{h}_i^{n+1} - h_i^n}{\Delta t} + \Ord{2} \\
     &= \frac{k}{M}\tilde{h}_i^{n+1} + \left( 1 - \frac{k}{M} \right) h_i^n + \Ord{2} \,.
\end{align*}
\end{linenomath}
Therefore, a second order approximation to \( \h \) data at times \( t^{n,k} \) for \( i \in \C_P^\ifone \) is given by
\begin{linenomath}
\begin{equation}
    \label{eqn:h_old_pred}
    h_i^{n,k} = \frac{k}{M} \tilde{h}_i^{n+1} + \left(1-\frac{k}{M}\right) h_i^n \,.
\end{equation}
\end{linenomath}
Similarly for \( \u^{n,k} \), take a Taylor series expansion centered at \( t^n \) for \( e \in \C_E^\ifone \),
\begin{linenomath}
\begin{align*}
    u_e\left( t^n + k\frac{\Delta t}{M} \right) &= u_e(t^n) + \left( t^n + k\frac{\Delta t}{M} - t^n \right) \pd{u_e}{t}(t^n) + \Ord{2} \\
    &=  u_e(t^n) + k\frac{\Delta t}{M} \pd{u_e}{t}(t^n) + \Ord{2} \,.
\end{align*}
\end{linenomath}
Approximate \( u_e(t^n) \) by \( u_e^n \), which is a second order approximation (the order of FB-RK(3,2)).
We approximate \(  \pd{u_e}{t}(t^n) \) similarly to the above, \( u_t(t^n) = \frac{\tilde{u}^{n+1} - u^n}{\Delta t} + \Ord{1} \).
Insert these approximations into the above to get 
\begin{linenomath}
\begin{align*}
    u_e\left( t^n + k\frac{\Delta t}{M} \right) &=  \left( u_e^n + \Ord{2} \right) + k\frac{\Delta t}{M} \left( \frac{\tilde{u}_e^{n+1} - u_e^n}{\Delta t} + \Ord{1} \right) + \Ord{2} \\
    &= u_e^n +  k\frac{\Delta t}{M} \frac{\tilde{u}_e^{n+1} - u_e^n}{\Delta t} + \Ord{2} \\
    &= \frac{k}{M}\tilde{u}_e^{n+1} + \left( 1 - \frac{k}{M} \right) u_e^n + \Ord{2} \,.
\end{align*}
\end{linenomath}
A second order approximation to \( \u \) data times \( t^{n,k}\) for \( e \in \C_E^\ifone \) is given by
\begin{linenomath}
\begin{equation}
    \label{eqn:u_old_pred}
    u_e^{n,k} = \frac{k}{M} \tilde{u}_e^{n+1} + \left(1-\frac{k}{M}\right) u_e^n \,.
\end{equation}
\end{linenomath}

Next, we require predictions for \( \bar{\h}^{n,k+\nicefrac{1}{3}} \) and \( \bar{\u}^{n,k+\nicefrac{1}{3}} \) for \( k = 0, \cdots, M-1 \).
Starting with \( \bar{\h}^{n,k+\nicefrac{1}{3}} \), we insert the above prediction \( \h^{n,k} \), given by \eqref{eqn:h_old_pred}, into the first stage of FB-RK(3,2) with the fine time-step \( \frac{\Delta t}{M} \) for \( i \in \C_P^\ifone \),
\begin{linenomath}
\begin{equation*}
    \bar{h}_i^{n,k+\nicefrac{1}{3}} = h_i^{n,k} + \frac{\Delta t}{3M} \Psi_i\left( \u^{n,k}, \h^{n,k} \right) \,.
\end{equation*}
\end{linenomath}
We need a way to approximate \(  \Psi_i\left( u^{n,k}, h^{n,k} \right) \).
Assuming sufficient smoothness of \( \Psi_i \), we can use that \(  \Psi_i\left( u^{n,k}, h^{n,k} \right) =  \Psi_i\left( u^n, h^n \right) + \Ord{1} \).
Insert this into the above and we get
\begin{linenomath}
\begin{align*}
    \bar{h}_i^{n,k+\nicefrac{1}{3}} &= h_i^{n,k} + \frac{\Delta t}{3M} \left( \Psi_i\left( u^n, h^n \right) + \Ord{1} \right) \\
    &= \left( \frac{k}{M}\tilde{h}_i^{n+1} + \left( 1 - \frac{k}{M} \right) h_i^n + \Ord{2} \right) + \frac{\Delta t}{3M} \Psi_i\left( \u^n, \h^n \right) + \Ord{2} \\
    &= \frac{k}{M}\tilde{h}_i^{n+1} + \left( 1 - \frac{k}{M} \right) h_i^n + \frac{\Delta t}{3M} \Psi_i\left( u^n, h^n \right) + \Ord{2} \,.
\end{align*}
\end{linenomath}
Then, we know that \(  \Psi_i\left( \u^n, \h^n \right) = \frac{3(\tilde{h}_i^{n+\nicefrac{1}{3}} - h_i^n)}{\Delta t} \).
Substitute this to get
\begin{linenomath}
\begin{align*}
     \bar{h}_i^{n,k+\nicefrac{1}{3}} &= \frac{k}{M}\tilde{h}_i^{n+1} + \left( 1 - \frac{k}{M} \right) h_i^n + \frac{\Delta t}{3M} \frac{3(\tilde{h}^{n+\nicefrac{1}{3}} - h^n)}{\Delta t} + \Ord{2} \\
     &= \frac{k}{M}\tilde{h}_i^{n+1} + \frac{1}{M} \tilde{h}_i^{n+\nicefrac{1}{3}} +  \left( 1 - \frac{k+1}{M} \right) h_i^n + \Ord{2} \,.
\end{align*}
\end{linenomath}
A second order approximation to first-stage \( \h \) data at intermediate time levels for \( i \in \C_P^\text{IF-1} \) is given by
\begin{linenomath}
\begin{equation}
    \label{eqn:h_s1_pred}
    \bar{h}_i^{n,k+\nicefrac{1}{3}} = \frac{k}{M}\tilde{h}_i^{n+1} + \frac{1}{M} \tilde{h}_i^{n+\nicefrac{1}{3}} +  \left( 1 - \frac{k+1}{M} \right) h_i^n \,.
\end{equation}
\end{linenomath}
Proceed similarly for \( \bar{\u}^{n,k+\nicefrac{1}{3}} \), insert \( \u^{n,k} \) into the first stage of FB-RK(3,2) with the fine time-step for \( e \in \C_E^\ifone \),
\begin{linenomath}
\begin{equation*}
    \bar{u}_e^{n,k+\nicefrac{1}{3}} = u_e^{n,k} + \frac{\Delta t}{3M} \Phi_e\left( \u^{n,k}, \h^{*,k} \right) \,.
\end{equation*}
\end{linenomath}
We need to approximate \(  \Phi_e\left( \u^{n,k}, \h^{*,k} \right) \).
We would like to do this by \( \Phi_e\left( \u^n, \h^* \right) = \frac{3(\tilde{u}_e^{n+\nicefrac{1}{3}} - u_e^n)}{\Delta t} \).
One can show that \( \abs{\Phi_e\left( \u^{n,k}, \h^{*,k}\right) - \Phi_e\left( \u^n, \h^* \right)} \leq \Ord{1} \), and so \( \Phi_e\left( \u^{n,k}, \h^{*,k}\right) = \Phi_e\left( \u^n, \h^* \right) + \Ord{1} \).
This follows from the assumption that \( \Phi_e \) is Lipschitz, and we omit the details for reasons of space.
Insert this approximation into the above and we get
\begin{linenomath}
\begin{align*}
    \bar{u}_e^{n,k+\nicefrac{1}{3}} &=  u_e^{n,k} + \frac{\Delta t}{3M}\left( \Phi_e\left( \u^n, \h^* \right) + \Ord{1} \right) \\
    &= \left( \frac{k}{M} \tilde{u}^{n+1} + \left(1-\frac{k}{M}\right) u^n + \Ord{2} \right) + \frac{\Delta t}{3M} \frac{3(\tilde{u}_e^{n+\nicefrac{1}{3}} - u_e^n)}{\Delta t} + \Ord{2} \\
    &= \frac{k}{M}\tilde{u}_e^{n+1} + \frac{1}{M} \tilde{u}_e^{n+\nicefrac{1}{3}} +  \left( 1 - \frac{k+1}{M} \right) u_e^n + \Ord{2} \,.
\end{align*}
\end{linenomath}
A second order approximation to first-stage \( \u \) data at intermediate time levels for \( e \in \C^\text{IF-1} \) is given by
\begin{linenomath}
\begin{equation}
    \label{eqn:u_s1_pred}
    \bar{u}_e^{n,k+\nicefrac{1}{3}} = \frac{k}{M}\tilde{u}_e^{n+1} + \frac{1}{M} \tilde{u}_e^{n+\nicefrac{1}{3}} +  \left( 1 - \frac{k+1}{M} \right) u_e^n \,.
\end{equation}
\end{linenomath}

Finally, we predictions for second-stage data \( \bar{\h}^{n,k+\nicefrac{1}{2}} \) and \( \bar{\u}^{n,k+\nicefrac{1}{2}} \) for \( k = 0, \cdots, M-1 \).
Starting with \( \bar{\h}^{n,k+\nicefrac{1}{3}} \), we insert the above prediction \( \h^{n,k} \) into the second stage of FB-RK(3,2) with the fine time-step \( \frac{\Delta t}{M} \) for \( i \in \C_P^\ifone \),
\begin{linenomath}
\begin{equation*}
    \bar{h}_i^{n,k+\nicefrac{1}{2}} = h_i^{n,k} + \frac{\Delta t}{2M} \Psi_i\left( \bar{\u}^{n,k+\nicefrac{1}{3}}, \bar{\h}^{n,k+\nicefrac{1}{3}} \right) \,.
\end{equation*}
\end{linenomath}
We need a way to approximate \(  \Psi_i\left( \bar{\u}^{n,k+\nicefrac{1}{3}}, \bar{\h}^{n,k+\nicefrac{1}{3}} \right) \).
Again using the assumption that we have Lipschitz continuity, we can use that \(  \Psi_i\left( \bar{\u}^{n,k+\nicefrac{1}{3}}, \bar{\h}^{n,k+\nicefrac{1}{3}} \right) =  \Psi_i\left( \tilde{\u}^{n+\nicefrac{1}{3}}, \tilde{\h}^{n+\nicefrac{1}{3}} \right) + \Ord{1} \).
Insert this into the above and we get
\begin{linenomath}
\begin{align*}
    \bar{h}_i^{n,k+\nicefrac{1}{2}} &= h_i^{n,k} + \frac{\Delta t}{2M} \left( \Psi_i\left(  \tilde{\u}^{n+\nicefrac{1}{3}}, \tilde{\h}^{n+\nicefrac{1}{3}} \right) + \Ord{1} \right) \\
    &= \left( \frac{k}{M}\tilde{h}_i^{n+1} + \left( 1 - \frac{k}{M} \right) h_i^n + \Ord{2} \right) +  \frac{\Delta t}{2M} \frac{2(\tilde{h}_i^{n+\nicefrac{1}{2}} - h_i^n)}{\Delta t} + \Ord{2} \\
    &= \frac{k}{M}\tilde{h}_i^{n+1} + \frac{1}{M} \tilde{h}_i^{n+\nicefrac{1}{2}} +  \left( 1 - \frac{k+1}{M} \right) h_i^n + \Ord{2} \,.
\end{align*}
\end{linenomath}
A second-order approximation to stage-two \( \h \) data for \( i \in \C_P^\ifone \) is given by
\begin{linenomath}
\begin{equation}
    \label{eqn:h_s2_pred}
    \bar{h}_i^{n,k+\nicefrac{1}{2}} = \frac{k}{M}\tilde{h}_i^{n+1} + \frac{1}{M} \tilde{h}_i^{n+\nicefrac{1}{2}} +  \left( 1 - \frac{k+1}{M} \right) h_i^n \,.
\end{equation}
\end{linenomath}
For \( \bar{u}_e^{n,k\nicefrac{1}{2}} \) for \( e \in \C_E^\ifone \), proceed similarly:
\begin{linenomath}
\begin{equation*}
    \bar{u}_e^{n,k+\nicefrac{1}{2}} = u_e^{n,k} + \frac{\Delta t}{2M} \Phi_e\left( \bar{\u}^{n,k+\nicefrac{1}{3}}, \h^{**,k} \right) \,.
\end{equation*}
\end{linenomath}
As we did to obtain the prediction for \( \bar{h}_i^{n,k+\nicefrac{1}{2}} \), we can approximate \( \Phi_e\left( \bar{\u}^{n,k+\nicefrac{1}{3}}, \h^{**,k} \right) = \Phi_e\left( \tilde{\u}^{n+\nicefrac{1}{3}}, \h^{**} \right) + \Ord{1} \). 
Insert this into the above to get
\begin{linenomath}
\begin{align*}
     \bar{u}_e^{n,k+\nicefrac{1}{2}} &= u_e^{n,k} + \frac{\Delta t}{2M} \left( \Phi_e\left( \tilde{\u}^{n+\nicefrac{1}{3}}, \h^{**} \right) + \Ord{1} \right) \\
     &= \left(  \frac{k}{M} \tilde{u}^{n+1} + \left(1-\frac{k}{M}\right) u^n + \Ord{2} \right) + \frac{\Delta t}{2M} \frac{2(\tilde{u}_e^{n+\nicefrac{1}{2}} - u_e^n)}{\Delta t} + \Ord{2} \\
     &= \frac{k}{M}\tilde{u}_e^{n+1} + \frac{1}{M} \tilde{u}_e^{n+\nicefrac{1}{2}} +  \left( 1 - \frac{k+1}{M} \right) u_e^n + \Ord{2} \,.
\end{align*}
\end{linenomath}
A second-order approximation to stage-two \( \u \) data for \( e \in \C_E^\ifone \) is given by
\begin{linenomath}
\begin{equation}
    \label{eqn:u_s2_pred}
    \bar{u}_e^{n,k+\nicefrac{1}{2}} = \frac{k}{M}\tilde{u}_e^{n+1} + \frac{1}{M} \tilde{u}_e^{n+\nicefrac{1}{2}} +  \left( 1 - \frac{k+1}{M} \right) u_e^n \,.
\end{equation}
\end{linenomath}

Taking all this together, the interface one prediction step is given by \eqref{eqn:preds}.


\end{document}



\bibliographystyle{elsarticle-harv}
\bibliography{references}


\end{document}